\documentclass[a4paper,11pt,reqno]{amsart}
\usepackage{pdfsync}

\def\usealternativegraphics{false}
\def\usetodonotes{false}

\usepackage{setspace}

\usepackage{color}

\usepackage{verbatim,ifthen}
\usepackage{amsmath,amssymb,latexsym,ifpdf}
\usepackage{pgf}

\ifthenelse{\equal{\usealternativegraphics}{false}}{\usepackage{tikz} 
\usetikzlibrary{calc}}{}

\newcounter{number_of_todos}
\ifthenelse{\equal{\usetodonotes}{true}}{
\usepackage{todonotes}

}{

}

\ifthenelse{\equal{\usealternativegraphics}{true}}{
\newenvironment{alternativegraphic}[2][1]{
\ifthenelse{\equal{none}{#2}}{
  \textbf{Picture omitted}
}{
  \ifpdf
    \IfFileExists{#2.pdf}{\includegraphics[page=#1]{#2.pdf}}{\textbf{PDF file not found.}}
  \else
    \IfFileExists{#2.eps}{\includegraphics{#2.eps}}{
      \IfFileExists{#2_#1.eps}{\includegraphics{#2_#1.eps}}{\textbf{EPS file not found.}}
    }
  \fi
}
\comment
}{
\endcomment
}
}{
\newenvironment{alternativegraphic}[2][]{}{}
}

\newtheorem{theorem}{Theorem}[section]
\newtheorem{lemma}[theorem]{Lemma}
\newtheorem{prop}[theorem]{Proposition}
\newtheorem{cor}[theorem]{Corollary}

\newtheorem*{prop*}{Proposition}
\newtheorem*{lemma*}{Lemma}
\theoremstyle{remark}
\newtheorem{definition}[theorem]{Definition}

\newtheorem{example}[theorem]{Example}
\newtheorem{remark}[theorem]{Remark}

\numberwithin{equation}{section}

\def\tr{\mathrm{Tr}}
\def\supp{\mathrm{supp}}
\def\ind{\mathrm{Ind}}
\def\L{\mathrm{L}}
\def\U{\mathrm{U}}
\def\M{\mathrm{M}}
\def\liminf{\mathrm{lim\, inf}}

\def\CC{\mathbb C}
\def\NN{\mathbb N}
\def\PP{\mathbb P}

\def\ZZ{\mathbb Z}
\def\RR{\mathbb R}

\def\ge{\geqslant}
\def\le{\leqslant}

\def\Nn{\mathcal N}
\def\Hh{\mathcal H}

\newcommand{\Ind}{\operatorname{Ind}}

\providecommand{\abs}[1]{\lvert#1\rvert}

\providecommand{\braces}[1]{\{#1\}}

\providecommand{\bkt}[1]{(#1)}

\title{Strength of convergence in the orbit space of a groupoid}
\author{Robert Hazlewood} 
\address{School of Mathematics and Statistics\\
The University of New South Wales\\Sydney\\
NSW 2052\\
Australia}
\email{robbiehazlewood@gmail.com}

\author{Astrid an Huef}
\address{Department of Mathematics and Statistics\\
University of Otago\\PO Box 56\\ Dunedin 9054\\
New Zealand}
\email{astrid@maths.otago.ac.nz}

\date{16 June 2010}

\begin{document}
\begin{abstract}
Let $G$ be a second-countable locally-compact Hausdorff groupoid with a Haar system, and let $\{x_n\}$ be a sequence in the unit space $G^{(0)}$ of $G$.  We show that the notions of strength of convergence of  $\{x_n\}$ in the orbit space  $G^{(0)}/G$ and  measure-theoretic accumulation along the orbits are equivalent ways of realising  multiplicity numbers associated to a sequence of induced  representation of the groupoid $C^*$-algebra. 
\end{abstract}

\thanks{This research was supported by the Australian Research Council.  Robert Hazlewood thanks the Department of Mathematics and Statistics at the University of Otago for their kind hospitality.}

\maketitle


\section{Introduction}\label{sec_intro}
Suppose $H$ is a locally-compact Hausdorff group acting freely and continuously on a  locally-compact Hausdorff space $X$, so that $(H,X)$ is a free transformation group. In \cite[pp.~95--96]{Green1977} Green gives an example of a free non-proper action of $H=\RR$ on a subset $X$ of $\RR^3$; the non-properness comes down to the existence of $z\in X$, $\braces{x_n}\subset X$, and two sequences $\braces{s_n}$ and $\braces{t_n}$ in $H$ such that
\begin{enumerate}\renewcommand{\theenumi}{\roman{enumi}}\renewcommand{\labelenumi}{(\theenumi)}
\item $s_n^{-1}\cdot x_n\rightarrow z$ and $t_n^{-1}\cdot x_n\rightarrow z$; and 
\item $t_ns_n^{-1}\rightarrow\infty$ as $n\rightarrow\infty$, in the sense that $\braces{t_ns_n^{-1}}$ has no convergent subsequence.
\end{enumerate}
In  \cite[Definition~2.2]{Archbold-Deicke2005}), and subsequently in \cite[p.~2]{Archbold-anHuef2006}, the sequence $\{x_n\}$ is said to converge $2$-times in the orbit space to $z\in X$.  Each orbit $H\cdot x$  gives an induced  representation $\Ind\epsilon_{x}$ of the associated transformation-group $C^*$-algebra $C_0(X)\rtimes H$ which is irreducible, and the $k$-times convergence of $\{x_n\}$ in the orbit space to $z\in X$ translates into  statements about  various multiplicity numbers associated to $\Ind\epsilon_z$ in the spectrum of $C_0(X)\rtimes H$, as in  \cite[Theorem~2.5]{Archbold-Deicke2005}, \cite[Theorem~1.1]{Archbold-anHuef2006} and \cite[Theorem~2.1]{Archbold-anHuef2008}.

Upper and lower multiplicity numbers associated to irreducible representations $\pi$ of a $C^*$-algebra $A$ were introduced by Archbold \cite{Archbold1994} and extended to multiplicity numbers relative to a net of irreducible representations  by Archbold and Spielberg  \cite{Archbold-Spielberg1996}. 
The upper multiplicity $M_U(\pi)$ of $\pi$, for example, counts  `the number of nets of orthogonal equivalent pure states which can converge to a common pure state associated to $\pi$' \cite[p. 26]{aklss2001}.  The definition of $k$-times convergence and \cite[Theorem~2.5]{Archbold-Deicke2005} were very much motivated by a notion of $k$-times convergence in the dual space of a nilpotent Lie group \cite{Ludwig1990} and its connection with relative multiplicity numbers (see, for example,  \cite[Theorem~2.4]{aklss2001} and  \cite[Theorem~5.8]{Archbold-Ludwig-Schlichting2007}).

Theorem~1.1 of \cite{Archbold-anHuef2006} shows that the topological property of a  sequence $\{x_n\}$ converging $k$-times in the orbit space to $z\in X$ is equivalent to (1) a measure theoretic accumulation along the orbits $G\cdot x_n$ and (2) that the lower  multiplicity of $\Ind\epsilon_z$ relative to the sequence $\{\Ind\epsilon_{x_n}\}$ is at least $k$. In this paper we prove that the results of \cite{Archbold-anHuef2006} generalise to principal groupoids.

In our main arguments we have tried to preserve as much as possible the structure of those in  \cite{Archbold-anHuef2006}, although the arguments presented here are often more complicated in order to cope with the partially defined product in a groupoid and the set of measures that is a Haar system compared to the fixed Haar measure  used in the transformation-group case. 

Our theorems have led us to a new class of examples exhibiting $k$-times convergence in groupoids that are not based on transformation groups, thus justifying our level of generality. Given a row-finite directed graph $E$, Kumjian, Pask, Raeburn and Renault in \cite{kprr1997} used the set of all infinite paths in $E$ to construct an r-discrete groupoid $G_E$, called a {\em path groupoid}.  We prove that $G_E$ is principal if and only if $E$ contains no cycles (Proposition~\ref{prop_principal_iff_no_cycles}). We then exhibit principal $G_E$ with Hausdorff and non-Hausdorff orbits space, respectively, both with a  $k$-times converging sequence  in the orbit space. In particular, our examples can be used to find a groupoid $G_E$ whose  $C^*$-algebra has non-Hausdorff spectrum and distinct upper and lower multiplicity counts among its irreducible representations.

\section{Preliminaries}\label{sec_prelim}

We denote the unit space of a groupoid $G$ by $G^{(0)}$. For $x\in G^{(0)}$ we call the set $r\big(s^{-1}(\{x\})\big)=s\big(r^{-1}(\{x\})\big)$ the {\em orbit} of $x$ and denote it by $[x]$. For a subset $U$ of $G^{(0)}$ we define $G_U:=s^{-1}(U)$, $G^U:=r^{-1}(U)$, and $G|_U:=s^{-1}(U)\cap r^{-1}(U)$. We denote the set of all positive integers by $\PP$ and the set of all non-negative integers by $\NN$.

\begin{definition}\label{def_Haar_system}A {\em right Haar system} on $G$ is a set $\{\lambda_x:x\in G^{(0)}\}$ of non-negative Radon measures on $G$ such that
\begin{enumerate}\renewcommand{\theenumi}{\roman{enumi}}
 \item $\supp\,\lambda_x=G_x$ $\big(=s^{-1}(\{x\})\big)$\quad for all $x\in G^{(0)}$;
\item\label{Haar_system_property_2} for $f\in C_c(G)$, the function $x\mapsto \int f\,d\lambda_x$ on $G^{(0)}$ is in $C_c(G^{(0)})$; and
\item\label{Haar_system_property_3} for $f\in C_c(G)$ and $\gamma\in G$, 
\[\int f(\alpha\gamma)\,d\lambda_{r(\gamma)}(\alpha)=\int f(\alpha)\,d\lambda_{s(\gamma)}(\alpha).\]
\end{enumerate}
We will refer to \eqref{Haar_system_property_2} as the {\em continuity of the Haar system} and to \eqref{Haar_system_property_3} as {\em Haar-system invariance}. The collection $\{\lambda^x:x\in G^{(0)}\}$ of measures where $\lambda^x(E):=\lambda_x(E^{-1})$ is a {\em left Haar system}, which is a system of measures such that $\supp\,\lambda^x=G^x$ and, for $f\in C_c(G)$, $x\mapsto \int f\,d\lambda^x$ is continuous and $\int f(\gamma\alpha)\,d\lambda^{s(\gamma)}(\alpha)=\int f(\alpha)\,d\lambda^{r(\gamma)}(\alpha)$. Given that we can easily convert a right Haar system $\{\lambda_x\}$ into a left Haar system $\{\lambda^x\}$ and vice versa, we will simply refer to a {\em Haar system} $\lambda$ and use subscripts to refer to elements of the right Haar system $\{\lambda_x\}$ and superscripts to refer to elements of the left Haar system $\{\lambda^x\}$.
\end{definition}

The following lemma follows from the invariance of the Haar system and the Dominated Convergence Theorem; we omit the proof.

\begin{lemma}[Haar-system invariance]
Suppose $G$ is a locally-compact Hausdorff groupoid with Haar system $\lambda$. If $K\subset G$ is compact and $\gamma\in G$ with $s(\gamma)=x$ and $r(\gamma)=y$, then $\lambda_x(K\gamma)=\lambda_y(K)$ and $\lambda^x(\gamma^{-1} K)=\lambda^y(K)$.
\end{lemma}

Definition~\ref{def_induced_representation} below is Definition 2.45 in the unpublished book \cite{Muhly-book}. Alternative descriptions of the induced representation may be found in \cite[pp.~234]{Muhly-Williams1990} and \cite[pp.~81--82]{Renault1980}.
\begin{definition}\label{def_induced_representation}
Suppose $G$ is a second-countable locally-compact Hausdorff groupoid with Haar system $\lambda$ and let $\mu$ be a Radon measure on $G^{(0)}$.
\begin{enumerate}\renewcommand{\labelenumi}{(\roman{enumi})}
\item We write $\nu=\mu\circ\lambda=\int \lambda^x \, d\mu$ for the measure on $G$ defined for every Borel-measurable function $f:G\rightarrow\CC$ by $\int_G f(\gamma)\, d\nu(\gamma)=\int_{G^{(0)}}\int_G f(\gamma)\, d\lambda^x(\gamma)\, d\mu(x)$. We call $\nu$ the measure induced by $\mu$, and we write $\nu^{-1}$ for the image of $\nu$ under the homeomorphism $\gamma\mapsto \gamma^{-1}$.
\item For $f\in C_c(G)$, $\ind\,\mu(f)$ is the operator on $L^2(G,\nu^{-1})$ defined by the formula
\begin{align*}
\big(\ind\,\mu(f)\xi\big)(\gamma)&=\int_G f(\alpha)\xi(\alpha^{-1}\gamma)\, d\lambda^{r(\gamma)}(\alpha)\\
&=\int_G f(\gamma \alpha)\xi(\alpha^{-1})\, d\lambda^{s(\gamma)}(\alpha).
\end{align*}
\end{enumerate}
\end{definition}
In this paper we are interested in representations that are induced by point-mass measures $\delta_x$ on $G^{(0)}$. We denote $\ind\,\delta_x$ by $\L^x$ for all $x\in G^{(0)}$ as in \cite{Muhly-Williams1990} and \cite{Clark-anHuef2008}. 

\begin{remark}\label{measure_induced_epsilon_x} It follows from the definition of the induced mesasure that for $x\in G^{(0)}$, the measure $\nu$ induced by $\delta_x$ is equal to $\lambda^x$. In particular we have $\nu^{-1}=\lambda_x$, so $\L^x$ acts on $L^2(G,\lambda_x)$. The operator $\L^x$ is then given by
\[
\big(\L^x(f)\xi\big)(\gamma)=\int_G f(\gamma \alpha^{-1})\xi(\alpha)\, d\lambda_x(\alpha)
\]
for all $\xi\in L^2(G,\lambda_x)$ and all $\gamma\in G$. There is a close relationship between the convolution on $C_c(G)$ and these induced representations: recall that for $f,g\in C_c(G)$, the convolution $f\ast g\in C_c(G)$ is given by
\[
 f\ast g(\gamma)=\int_G f(\gamma\alpha^{-1})g(\alpha)\,d\lambda_{s(\gamma)}(\alpha)\quad\text{for all }\gamma\in G,
\]
so that
\[
\big(\L^x(f)g\big)(\gamma)=f\ast g(\gamma)\quad\text{for any }x\in G^{(0)}\text{ and }\gamma\in G_x.
\]
We denote the norm in $L^2(G,\lambda_x)$ by $\|\cdot\|_x$. Finally note that when $G$ is a second-countable locally-compact principal groupoid that admits a Haar system, each $\L^x$ is irreducible by \cite[Lemma~2.4]{Muhly-Williams1990}.
\end{remark}
\begin{remark} If
 $G=(H, X)$ is a second-countable free transformation group, then the representations $\L^x$ defined above are unitarily equivalent to the representations $\Ind\epsilon_x$ used in \cite{Archbold-anHuef2006}. Specifically,   let $\nu$ be a choice of right Haar measure on $H$ and $\Delta$  the associated modular function. The map $\iota:C_c(H\times X)\to C_c(H\times X)$ defined by
\[
\iota(f)(t,x)=f(t,x)\Delta(t)^{1/2}
\]
extends to an isomorphism $\iota$ of the groupoid $C^*$-algebra $C^*(H\times X)$ onto the transformation-group $C^*$-algebra $C_0(X)\rtimes H$ \cite[p.~58]{Renault1980}. Fix $x\in X$.  Then there is a unitary  $U_x:L^2(H,\nu)\to L^2(H\times X,\lambda_x)$,  characterised by $U(\xi)(h, y)=\xi(h)\delta_x(h^{-1}\cdot y)$ for $\xi\in C_c(H)$, and 
$U(\Ind\epsilon_x(\iota(f))U^*=\L^x(f)$ for $f\in C^*(H\times X)$.
\end{remark}

Let $A$ be a $C^*$-algebra. We  write $\theta$ for the canonical surjection from the space $P(A)$ of pure states of $A$ to the spectrum $\hat A$ of $A$.  We  frequently  identify an irreducible representation $\pi$ with its equivalence class in $\hat A$ and we write $\Hh_\pi$ for the Hilbert space on which $\pi(A)$ acts.

Let $\pi\in\hat A$ and let $\braces{\pi_\alpha}$ be a net in $\hat{A}$. We now recall the definitions of {\em upper} and {\em lower multiplicity}  $M_\U(\pi)$ and $M_\L(\pi)$ from \cite{Archbold1994}, and {\em relative upper} and {\em relative lower multiplicity} $\M_\U(\pi,\braces{\pi_\alpha})$ and  $\M_\L(\pi,\braces{\pi_\alpha})$ from \cite{Archbold-Spielberg1996}. 

Let $\Nn$ be the weak$^*$-neighbourhood base at zero in the dual $A^*$ of $A$ consisting of all open sets of the form
\[
N=\{\psi\in A^*:|\psi(a_i)|<\epsilon, 1\le i\le n\},
\]
where $\epsilon>0$ and $a_1,a_2,\ldots,a_n\in A$.
Suppose $\phi$ is a pure state of $A$ associated with $\pi$ and let $N\in\Nn$. Define
\[
V(\phi,N)=\theta\big((\phi+N)\cap P(A)\big),
\]
an open neighbourhood of $\pi$ in $\hat{A}$. For $\sigma\in \hat{A}$ let
\[
\mathrm{Vec}(\sigma,\phi,N)=\braces{\eta\in \Hh_\sigma:\|\eta\|=1,(\sigma(\cdot)\eta\, |\,\eta)\in \phi+N}.
\]
Note that $\mathrm{Vec}(\sigma,\phi,N)$ is non-empty if and only if $\sigma\in V(\phi,N)$. For any $\sigma\in V(\phi,N)$  define $d(\sigma,\phi,N)$ to be the supremum in $\PP\cup\braces{\infty}$ of the cardinalities of finite orthonormal subsets of $\mathrm{Vec}(\sigma,\phi,N)$. Write $d(\sigma,\phi,N)=0$ when $\mathrm{Vec}(\sigma,\phi,N)$ is empty.

Define
\[
\M_\U(\phi,N)=\underset{\sigma\in V(\phi,N)}{\sup} d(\sigma,\phi,N)\in\PP\cup\braces{\infty}.
\]
Note that if $N'\in \Nn$ and $N\subset N$, then $\M_\U(\phi,N')\le\M_\U(\phi,N)$. Now define
\[
\M_\U(\phi)=\underset{N\in\Nn}{\inf}\M_\U(\phi,N)\in\PP\cup\braces{\infty}.
\]
By \cite[Lemma~2.1]{Archbold1994}, the value of $\M_\U(\phi)$ is independent of the pure state $\phi$ associated to $\pi$. Thus  $\M_\U(\pi):=\M_\U(\phi)$ is well-defined.  For lower multiplicity,  assume that $\braces{\pi}$ is not open in $\hat{A}$, and using \cite[Lemma~2.1]{Archbold1994} again,   define
\[
\M_\L(\pi):=\underset{N\in\Nn}\inf \Big(\underset{\sigma\rightarrow\pi, \sigma\ne\pi}{\lim\,\inf} d(\sigma,\phi,N)\Big)\in\PP\cup\braces{\infty}.
\]

Now suppose that $\{\pi_\alpha\}_{\alpha\in\Lambda}$ is a net in $\hat{A}$. For $N\in\Nn$ let
\[
\M_\U\big(\phi,N,\{\pi_\alpha\}\big)=\underset{\alpha\in\Lambda}{\lim\,\sup}\, d(\pi_\alpha,\phi,N)\in\NN\cup\braces{\infty}.
\]
Note that if $N'\in\Nn$ and $N'\subset N$ then $\M_\U\big(\phi,N',\{\pi_\alpha\}\big)\le \M_\U\big(\phi,N,\{\pi_\alpha\}\big)$. Then
\[
\M_\U\big(\pi,\{\pi_\alpha\}\big):=\underset{N\in\Nn}\inf \M_\U\big(\phi,N,\{\pi_\alpha\}\big)\in\NN\cup\braces{\infty},
\]
is well-defined because  the right-hand side  is independent of the choice of $\phi$ by an argument similar to the proof of \cite[Lemma~2.1]{Archbold1994}. Similarly define
\[
\M_\L\big(\phi,N,\{\pi_\alpha\}\big):=\underset{\alpha\in\Lambda}{\lim\,\inf}\, d(\pi_\alpha,\phi,N)\in\NN\cup\braces{\infty},
\]
and let
\[
\M_\L\big(\pi,\{\pi_\alpha\}\big)=\underset{N\in\Nn}\inf \M_\L\big(\phi,N,\{\pi_\alpha\}\big)\in\NN\cup\braces{\infty}.
\]
It follows that for any irreducible representation $\pi$ and any net $\{\pi_\alpha\}_{\alpha\in\Lambda}$ of irreducible representations,
\[
\M_\L\big(\pi,\{\pi_\alpha\}\big)\le\M_\U\big(\pi,\{\pi_\alpha\}\big)\le\M_\U(\pi)
\]
and, if $\{\pi_\alpha\}$ converges to $\pi$ with $\pi_\alpha\ne\pi$ eventually,
\[
\M_\L(\pi)\le\M_\L\big(\pi,\{\pi_\alpha\}\big).
\]
Finally, if $\{\pi_\beta\}$ is a subnet of $\{\pi_\alpha\}$, then
\[
\M_\L\big(\pi,\{\pi_\alpha\}\big)\le \M_\L\big(\pi,\{\pi_\beta\}\big)\le \M_\U\big(\pi,\{\pi_\beta\}\big)\le \M_\U\big(\pi,\{\pi_\alpha\}\big).
\]

\section{Lower multiplicity and $k$-times convergence I}\label{sec_lower_multiplicity_1}
A key goal for this paper is to describe the relationship between multiplicities of induced representations and strength of convergence in the orbit space. We start  this section by recalling the definition of $k$-times convergence in a groupoid from \cite{Clark-anHuef2008}. We then show that if a sequence converges $k$-times in the orbit space of a principal groupoid, then the lower multiplicity of the associated sequence of representations is at least $k$; the converse will be shown in Section \ref{sec_lower_multiplicity_2}.

Recall that a sequence $\{\gamma_n\}\subset G$ tends to infinity if it admits no convergent subsequence.
\begin{definition}\label{def_k-times_convergence}
Let $k\in\PP$. A sequence $\{x_n\}$ in $G^{(0)}$ is $k$-times convergent in $G^{(0)}/G$ to $z\in G^{(0)}$ if there exist $k$ sequences $\{\gamma_n^{(1)}\},\{\gamma_n^{(2)}\},\ldots,\{\gamma_n^{(k)}\}\subset G$ such that
\begin{enumerate}\renewcommand{\theenumi}{\roman{enumi}}\renewcommand{\labelenumi}{(\theenumi)}
\item\label{def_k-times_convergence_1} $s(\gamma_n^{(i)})=x_n$ for all $n$ and $1\le i\le k$;
\item\label{def_k-times_convergence_2} $r(\gamma_n^{(i)})\rightarrow z$ as $n\rightarrow\infty$ for $1\le i\le k$; and
\item\label{def_k-times_convergence_3} if $1\le i<j\le k$ then $\gamma_n^{(j)}(\gamma_n^{(i)})^{-1}\rightarrow\infty$ as $n\rightarrow\infty$.
\end{enumerate}
\end{definition}

The proof of the following proposition is based on \cite[Theorem~2.3]{Archbold-Deicke2005} and a part of \cite[Theorem~1.1]{Archbold-anHuef2006}.
\begin{prop}\label{AaH_thm_1.1_1_implies_2}
Suppose $G$ is a second-countable locally-compact Hausdorff principal groupoid with Haar system $\lambda$. Let $z\in G^{(0)}$ and suppose that $\{x_n\}$ is a sequence in $G^{(0)}$ that converges $k$-times to $z$ in $G^{(0)}/G$. Then $\M_\L(\L^z,\{\L^{x_n}\})\ge k$.
\begin{proof}
We will use a contradiction argument. Suppose that $\M_\L(\L^z,\{\L^{x_n}\})=r<k$. 
Fix a real-valued $g\in C_c(G)$ so that $\|g\|_z>0$. Define $\eta\in L^2(G,\lambda_z)$ by
\[
\eta(\alpha)=\|g\|_z^{-1}g(\alpha)\quad\text{for all }\alpha\in G. 
\]
Then
\[
\|\eta\|_z^2=\|g\|_z^{-2}\int g(\alpha)^2\, d\lambda_z(\alpha)=\|g\|_z^{-2}\|g\|_z^2=1,
\]
so $\eta$ is a unit vector in $L^2(G,\lambda_z)$
and the GNS construction of $\phi:=(\L^z(\cdot)\eta\,|\,\eta)$ is unitarily equivalent to $\L^z$. By the definition of lower multiplicity we now have
\[
\M_\L(\L^z,\{\L^{x_n}\})=\inf_{N\in\Nn} \M_\L(\phi,N,\{\L^{x_n}\})=r,
\]
so there exists $N\in\Nn$ such that
\[
\M_\L(\phi,N,\{\L^{x_n}\})=\underset{\scriptstyle n}{\lim\,\inf}\, d(\L^{x_n},\phi,N)=r,
\]
and consequently there exists a subsequence $\{y_m\}$ of $\{x_n\}$ such that
\begin{equation}\label{observation_to_contradict}
d(\L^{y_m},\phi,N)=r\quad\text{for all }m.
\end{equation}
 Since any subsequence of a sequence that is $k$-times convergent is also $k$-times convergent, we know that $\{y_m\}$ converges $k$-times to $z$ in $G^{(0)}/G$.

We will now use the $k$-times convergence of $\{y_m\}$ to construct $k$ sequences of unit vectors with sufficient properties to establish our contradiction. By the $k$-times convergence of $\{y_m\}$ there exist $k$ sequences
\[
\{\gamma_m^{(1)}\},\{\gamma_m^{(2)}\},\ldots,\{\gamma_m^{(k)}\}\subset G
\]
satisfying
{\allowdisplaybreaks\begin{enumerate}\renewcommand{\labelenumi}{(\roman{enumi})}
\item $s(\gamma_m^{(i)})=y_m$ for all $m$ and $1\le i\le k$;
\item $r(\gamma_m^{(i)})\rightarrow z$ as $m\rightarrow\infty$ for $1\le i\le k$; and
\item if $1\le i<j\le k$ then $\gamma_m^{(j)}(\gamma_m^{(i)})^{-1}\rightarrow\infty$ as $m\rightarrow\infty$.
\end{enumerate}}

For each $1\le i\le k$ and $m\ge 1$, define $\eta_m^{(i)}$ by
\[
\eta_m^{(i)}(\alpha)=\|g\|_{r(\gamma_m^{(i)})}^{-1}g\big(\alpha(\gamma_m^{(i)})^{-1}\big)\quad\text{for all }\alpha\in G.
\]
It follows from Haar-system invariance that
\begin{align*}
\|\eta_m^{(i)}\|_{y_m}^2&=\|g\|_{r(\gamma_m^{(i)})}^{-2}\int g\big(\alpha(\gamma_m^{(i)})^{-1}\big)^2\,d\lambda_{y_m}(\alpha)\\
&=\|g\|_{r(\gamma_m^{(i)})}^{-2}\int g(\alpha)^2\,d\lambda_{r(\gamma_m^{(i)})}(\alpha)\\
&=\|g\|_{r(\gamma_m^{(i)})}^{-2}\|g\|_{r(\gamma_m^{(i)})}^2=1,
\end{align*}
so the $\eta_m^{(i)}$ are unit vectors in $L^2(G,\lambda_{y_m})$. 

Now suppose that $1\le i<j\le k$. Then
\begin{equation}\label{working_the_etas}
(\eta_m^{(i)}\,|\,\eta_m^{(j)})_{y_m}= \|g\|_{r(\gamma_m^{(i)})}^{-1}\|g\|_{r(\gamma_m^{(j)})}^{-1}\int g\big(\alpha(\gamma_m^{(i)})^{-1}\big)g\big(\alpha(\gamma_m^{(j)})^{-1}\big)\, d\lambda_{y_m}(\alpha)\\
\end{equation}
Since $\gamma_m^{(i)}(\gamma_m^{(j)})^{-1}\rightarrow\infty$, $\gamma_m^{(i)}(\gamma_m^{(j)})^{-1}$ is eventually not in the compact set $(\supp\, g)(\supp\, g)^{-1}$, and so there exists $m_0$ such that if $m\ge m_0$, then
\[
(\supp\, g)\gamma_m^{(i)}\cap(\supp\, g)\gamma_m^{(j)}=\emptyset.
\]
(To see this claim, note that if $(\supp\, g)\gamma_m^{(i)}\cap (\supp\, g)\gamma_m^{(j)}\ne\emptyset$ then there exist $\alpha,\beta\in \supp\, g$ such that $\alpha\gamma_m^{(i)}=\beta\gamma_m^{(j)}$, and so $\gamma_m^{(i)}(\gamma_m^{(j)})^{-1}=\alpha^{-1}\beta\in (\supp\, g)^{-1}(\supp\, g)$.) For the integrand of \eqref{working_the_etas} to be non-zero, both $\alpha(\gamma_m^{(i)})^{-1}$ and $\alpha(\gamma_m^{(j)})^{-1}$ must be in $\supp\, g$, so $\alpha$ must be in $(\supp\, g)\gamma_m^{(i)}\cap(\supp\, g)\gamma_m^{(j)}$. But this is not possible if $m\ge m_0$. Thus, for any distinct $i,j$, we will eventually have $\eta_m^{(i)}\perp\eta_m^{(j)}$.

For the last main component of this proof we will establish that
\[
\big(\L^{y_m}(\cdot)\eta_m^{(i)}\,\big|\,\eta_m^{(i)}\big)\rightarrow \big(\L^z(\cdot)\eta\,\big|\,\eta\big)=\phi\quad\text{as }m\rightarrow\infty
\]
in the dual of $C^*(G)$ with the weak$^*$ topology for each $i$. Fix $f\in C_c(G)$. We have
{\allowdisplaybreaks\begin{align}
\big(\L^{z}(f)\eta\,\big|\,\eta\big)&=
\int_G\big(\L^z(f)\eta\big)(\alpha)\eta(\alpha)\, d\lambda_z(\alpha)\notag\\
&=\int_G\int_G f(\alpha\beta^{-1})\eta(\beta)\eta(\alpha)\,d\lambda_z(\beta)\,d\lambda_z(\alpha)\notag\\
\label{dealing_with_etas} &=\|g\|_z^{-2}\int_G\int_G f(\alpha\beta^{-1})g(\beta)g(\alpha)\,d\lambda_z(\beta)\,d\lambda_z(\alpha)
\end{align}
}Now fix $1\le i\le k$. By the invariance of the Haar system we have
{\allowdisplaybreaks\begin{align}
\big(\L^{y_m}&(f)\eta_m^{(i)}\,\big|\,\eta_m^{(i)}\big)=\int_G\int_G f(\alpha\beta^{-1})\eta_m^{(i)}(\beta)\eta_m^{(i)}(\alpha)\,d\lambda_{y_m}(\beta)\,d\lambda_{y_m}(\alpha)\notag\\
&=\|g\|_{r(\gamma_m^{(i)})}^{-2}\int_G\int_G f(\alpha\beta^{-1})g\big(\alpha(\gamma_m^{(i)})^{-1}\big)g\big(\beta(\gamma_m^{(i)})^{-1}\big)\, d\lambda_{y_m}(\beta)\,d\lambda_{y_m}(\alpha)\notag\\
&=\|g\|_{r(\gamma_m^{(i)})}^{-2}\int_G\int_G f(\alpha\beta^{-1})g(\alpha)g(\beta)\, d\lambda_{r(\gamma_m^{(i)})}(\beta)\, d\lambda_{r(\gamma_m^{(i)})}(\alpha)\notag\\
&=\|g\|_{r(\gamma_m^{(i)})}^{-2}\int_G f\ast g(\alpha) g(\alpha)\, d\lambda_{r(\gamma_m^{(i)})}(\alpha).\label{still_working_etas}
\end{align}
}We know that $r(\gamma_m^{(i)})\rightarrow z$ as $m\rightarrow\infty$ so, by the continuity of the Haar system, $\|g\|_{r(\gamma_m^{(i)})}\rightarrow \|g\|_z$ as $m\rightarrow\infty$. Since $f\ast g\in C_c(G)$ we can apply the continuity of the Haar system with \eqref{dealing_with_etas} and \eqref{still_working_etas} to see that
\begin{align*}
\big(\L^{y_m}(f)\eta_m^{(i)}\,\big|\,\eta_m^{(i)}\big)&=\|g\|_{r(\gamma_m^{(i)})}^{-2}\int_G f\ast g(\alpha) g(\alpha)\, d\lambda_{r(\gamma_m^{(i)})}(\alpha)\\
&\rightarrow \|g\|_{z}^{-2}\int_G f\ast g(\alpha) g(\alpha)\, d\lambda_{z}(\alpha)=\big(\L^{z}(f)\eta\,\big|\,\eta\big)
\end{align*}
as $m\rightarrow\infty$.

%

We have thus shown that, for each $i$,
\[
\big(\L^{y_m}(\cdot)\eta_m^{(i)}\,\big|\,\eta_m^{(i)}\big)\rightarrow \big(\L^z(\cdot)\eta\,\big|\,\eta\big)=\phi
\]
in the dual of $C^*(G)$ equipped with the weak$^*$ topology. Thus there exists $m_1$ such that for any $m\ge m_1$ and any $1\le i\le k$, the pure state $\big(\L^{y_m}(\cdot)\eta_m^{(i)}\,\big|\,\eta_m^{(i)}\big)$ is in $\phi+N$.

We have now established that every $\eta_m^{(i)}$ with $m\ge\max\{m_0,m_1\}$ is in $\mathrm{Vec}(\L^{y_m},\phi,N)$ with $\eta_m^{(i)}\perp\eta_m^{(j)}$ for $i\ne j$, so $d(\L^{y_m},\phi,N)\ge k$ for all $m\ge\max\{m_0,m_1\}$, contradicting our earlier observation \eqref{observation_to_contradict} that $d(\L^{y_m},\phi,N)=r<k$ for all $m$.
\end{proof}
\end{prop}

\section{Measure ratios and $k$-times convergence} In this section we show that lower bounds on measure ratios along orbits   gives strength of convergence in the orbit space. We begin by generalising \cite[Proposition~4.1]{Archbold-anHuef2006}.

\begin{prop}\label{AaH_prop4_1_1}
Let $G$ be a second-countable locally-compact Hausdorff principal groupoid with Haar system $\lambda$. Let $k\in\PP$ and $z\in G^{(0)}$ with $[z]$ locally closed in $G^{(0)}$. Assume that $\{x_n\}$ is a sequence in $G^{(0)}$ such that $[x_n]\rightarrow [z]$ uniquely in $G^{(0)}/G$. Suppose $\{W_m\}$ is a basic decreasing sequence of compact neighbourhoods of $z$ such that each $m$ satisfies 
\[
\underset{\scriptstyle n}{\lim\,\inf}\,\lambda_{x_n}(G^{W_m})>(k-1)\lambda_z(G^{W_m}).
\]
Then $\{x_n\}$ converges $k$-times in $G^{(0)}/G$ to $z$.
\begin{proof}
Let $\{K_m\}$ be an increasing sequence of compact subsets of $G$ such that $G=\bigcup_{m\ge 1}\mathrm{Int}\,K_m$. By the regularity of $\lambda_z$, for each $m\ge 1$ there exist $\delta_m>0$ and an open neighbourhood $U_m$ of $G_z^{W_m}$ such that
\begin{equation}\label{AaH4.1}
\underset{\scriptstyle n}{\lim\,\inf}\,\lambda_{x_n}(G^{W_m})>(k-1)\lambda_z(U_m)+\delta_m.
\end{equation}

We will construct, by induction, a strictly increasing sequence of positive integers $\{n_m\}$ such that, for all $n\ge n_m$,
\begin{align}
&\lambda_{x_n}(K_m\alpha\cap G^{W_m})<\lambda_z(U_m)+\delta_m/k\quad\text{for all }\alpha\in G_{x_n}^{W_m},\quad\text{and}\label{AaH4.2}\\
&\lambda_{x_n}(G^{W_m})>(k-1)\lambda_z(U_m)+\delta_m.\label{AaH4.3}
\end{align}

By applying Lemma \ref{the_unbroken_lemma} with $\delta=\lambda_z(U_1)-\lambda_z(G^{W_1})+\delta_1/k$ there exists $n_1$ such that $n\ge n_1$ implies
\[\lambda_{x_n}(K_1\alpha\cap G^{W_1})<\lambda_z(U_1)+\delta_1/k \quad\text{for all }\alpha\in G_{x_n}^{W_1},\] establishing \eqref{AaH4.2} for $m=1$. If necessary we can increase $n_1$ to ensure \eqref{AaH4.3} holds for $m=1$ by considering \eqref{AaH4.1}. Assuming that we have constructed $n_1<n_2<\cdots<n_{m-1}$, we apply Lemma \ref{the_unbroken_lemma} with $\delta=\lambda_z(U_m)-\lambda_z(G^{W_m})+\delta_m/k$ to obtain $n_m>n_{m-1}$ such that \eqref{AaH4.2} holds, and again, if necessary, increase $n_m$ to obtain \eqref{AaH4.3}.

If $n_1>1$ then, for each $1\le n<n_1$ and $1\le i\le k$, let $\gamma_n^{(i)}=x_n$. For each $n\ge n_1$ there is a unique $m$ such that $n_m\le n<n_{m+1}$. For every such $n$ and $m$ choose $\gamma_n^{(1)}\in G_{x_n}^{W_m}$ (which is always non-empty by \eqref{AaH4.3}). Using \eqref{AaH4.2} and \eqref{AaH4.3} we have
\begin{align*}
\lambda_{x_n}(G^{W_m}\backslash K_m\gamma_n^{(1)})&=\lambda_{x_n}(G^{W_m})-\lambda_{x_n}(G^{W_m}\cap K_m\gamma_n^{(1)})\\
&>\big((k-1)\lambda_z(U_m)+\delta_m\big)-\big(\lambda_z(U_m)+\delta_m/k\big)\\
&=(k-2)\lambda_z(U_m)+\frac{(k-1)}k\delta_m.
\end{align*}
So for each $n\ge n_1$ and its associated $m$ we can choose $\gamma_n^{(2)}\in G_{x_n}^{W_m}\backslash K_m\gamma_n^{(1)}$. We now have
\begin{align*}
\lambda_{x_n}&\big(G^{W_m}\backslash (K_m\gamma_n^{(1)}\cup K_m\gamma_n^{(2)})\big)\\
&=\lambda_{x_n}(G^{W_m}\backslash K_m\gamma_n^{(1)})-\lambda_{x_n}\big((G^{W_m}\backslash K_m\gamma_n^{(1)})\cap K_m\gamma_n^{(2)}\big)\\
&\ge\lambda_{x_n}(G^{W_m}\backslash K_m\gamma_n^{(1)})-\lambda_{x_n}(G^{W_m}\cap K_m\gamma_n^{(2)})\\
&>\Big((k-2)\lambda_z(U_m)+\frac{(k-1)}k\delta_m\Big)-\Big(\lambda_z(U_m)+\delta_m/k\Big)\\
&=(k-3)\lambda_z(U_m)+\frac{(k-2)}k\delta_m,
\end{align*}
enabling us to choose $\gamma_n^{(3)}\in G_{x_n}^{W_m}\backslash (K_m\gamma_n^{(1)}\cup K_m\gamma_n^{(2)})$. By continuing this process, for each $j=3,\ldots,k$ and each $n\ge n_1$ we have
\[
\lambda_{x_n}\Bigg(G^{W_m}\backslash \bigg(\bigcup_{i=1}^{j-1}K_m\gamma_n^{(i)}\bigg)\Bigg)>(k-j)\lambda_z(U_m)+\frac{(k-j-1)\delta_m}k,
\]
enabling us to choose
\begin{equation}\label{eqn_choosing_gammas}
\gamma_n^{(j)}\in G^{W_m}_{x_n}\backslash \bigg(\bigcup_{i=1}^{j-1}K_m\gamma_n^{(i)}\bigg).
\end{equation}
Note that for $n_m\le n<n_{m+1}$ we have $\gamma_n^{(j)}\notin K_m\gamma_n^{(i)}$ for $1\le i<j\le k$.

We will now establish that $x_n$ converges $k$-times to $z$ in $G^{(0)}/G$ by considering the $\gamma_n^{(i)}$. Note that $s(\gamma_n^{(i)})=x_n$ for all $n$ and $i$ by our choice of the $\gamma_n^{(i)}$. We will now establish that $r(\gamma_n^{(i)})\rightarrow z$ as $n\rightarrow\infty$ for $1\le i\le k$. To see this, fix $i$ and let $V$ be an open neighbourhood of $z$. Since $W_m\rightarrow \{z\}$ there exists $m_0$ such that $m\ge m_0$ implies $W_m\subset V$. For each $n\ge n_{m_0}$ there exists a $m\ge m_0$ such that $n_m\le n <n_{m+1}$, and so $r(\gamma_n^{(i)})\in W_m\subset V$.

Finally we claim that $\gamma_n^{(j)}(\gamma_n^{(i)})^{-1}\rightarrow\infty$ as $n\rightarrow\infty$ for $1\le i<j\le k$. Fix $i<j$ and let $K$ be a compact subset of $G$. There exists $m_0$ such that $K\subset K_m$ for all $m\ge m_0$. Then for $n\ge n_{m_0}$ there exists $m\ge m_0$ such that $n_m\le n<n_{m+1}$. By \eqref{eqn_choosing_gammas} we know
\begin{align*}
\gamma_n^{(j)}&\in G_{x_n}^{W_m}\backslash (K_m\gamma_n^{(i)})\\
&=\big(G_{x_n}^{W_m}(\gamma_n^{(i)})^{-1}\gamma_n^{(i)}\big)\backslash (K_m\gamma_n^{(i)})\\
&=\big((G_{x_n}^{W_m}(\gamma_n^{(i)})^{-1})\backslash K_m\big)\gamma_n^{(i)},
\end{align*}
and so $\gamma_n^{(j)}(\gamma_n^{(i)})^{-1}\in \big(G_{x_n}^{W_m}(\gamma_n^{(i)})^{-1}\big)\backslash K_m\subset G\backslash K_m\subset G\backslash K$, enabling us to conclude that $\{x_n\}$ converges $k$-times in $G^{(0)}/G$ to $z$.
\end{proof}
\end{prop}

In Proposition~\ref{part_of_AaH4.2_generalisation} below we prove a generalisation of a part of \cite[Proposition~4.2]{Archbold-anHuef2006}; to do this we need the following two lemmas.

\begin{lemma}\label{astrids_lim_sup}
Suppose $G$ is a second-countable groupoid with Haar system $\lambda$ and let $K$ be a compact subset of $G$. If $\{x_n\}\subset G^{(0)}$ is a sequence that converges to $z\in G^{(0)}$, then
\[
\underset{\scriptstyle n}{\lim\,\sup}\,\lambda_{x_n}(K)\le\lambda_z(K).
\]
\begin{proof}
Fix $\epsilon>0$. By the outer regularity of $\lambda_z$, there exists an open neighbourhood $U$ of $K$ such that
\[
\lambda_z(K)\le \lambda_z(U)<\lambda_z(K)+\epsilon/2.
\]
By Urysohn's Lemma there exists $f\in C_c(G)$ with $0\le f\le 1$ such that $f$ is identically one on $K$ and zero off $U$. In particular we have
\begin{equation}\label{measure_of_f_near_measure_of_K}
\lambda_z(K)\le\int f\,d\lambda_z<\lambda_z(K)+\epsilon/2.
\end{equation}

The continuity of the Haar system implies $\int f\,d\lambda_{x_n}\rightarrow \int f\,d\lambda_z$, so there exists $n_0$ such that $n\ge n_0$ implies
\[
\int f\,d\lambda_z-\epsilon/2 < \int f\,d\lambda_{x_n}<\int f\,d\lambda_z +\epsilon/2.
\]
By our choice of $f$ we have $\lambda_{x_n}(K)\le \int f\,d\lambda_{x_n}$, so
\[
\lambda_{x_n}(K)\le\int f\,d\lambda_{x_n}<\int f\,d\lambda_z+\epsilon/2.
\]
Combining this with \eqref{measure_of_f_near_measure_of_K} enables us to observe that for $n\ge n_0$, $\lambda_{x_n}(K)<\lambda_z(K)+\epsilon$, completing the proof.
\end{proof}
\end{lemma}
\begin{lemma}\label{corollary_astrid_lemma}
Suppose $G$ is a second-countable groupoid with Haar system $\lambda$ and let $K$ be a compact subset of $G$. For every $\epsilon>0$ and $z\in G^{(0)}$ there exists a neighbourhood $U$ of $z\in G^{(0)}$ such that $x\in U$ implies $\lambda_x(K)<\lambda_z(K)+\epsilon$.
\begin{proof}
Fix $\epsilon>0$ and $z\in G^{(0)}$. Let $\{U_n\}$ be a decreasing neighbourhood basis for $z$ in $G^{(0)}$. If our claim is false, then each $U_n$ contains an element $x_n$ such that $\lambda_{x_n}(K)\ge\lambda_z(K)+\epsilon$. But since each $x_n\in U_n$, $x_n\rightarrow z$, and so by Lemma \ref{astrids_lim_sup} there exists $n_0$ such that $n\ge n_0$ implies $\lambda_{x_n}(K)<\lambda_z(K)+\epsilon$, a contradiction.
\end{proof}
\end{lemma}

\begin{prop}\label{part_of_AaH4.2_generalisation}
Suppose $G$ is a second-countable locally-compact Hausdorff groupoid with Haar system $\lambda$. Suppose that $z\in G^{(0)}$ with $[z]$ locally closed in $G^{(0)}$ and suppose $\{x_n\}$ is a sequence in $G^{(0)}$. Assume that for every open neighbourhood $V$ of $z$ in $G^{(0)}$ such that $G_z^V$ is relatively compact, $\lambda_{x_n}(G^V)\rightarrow \infty$ as $n\rightarrow\infty$. Then, for every $k\ge 1$, the sequence $\{x_n\}$ converges $k$-times in $G^{(0)}/G$ to $z$.
\begin{proof}
Let $\{ K_m\}$ be an increasing sequence of compact subsets of $G$ such that $G=\cup_{m\ge 1}\,\mathrm{Int}\, K_m$. By Lemma \ref{corollary_astrid_lemma}, for each $K_m$ there exists an open neighbourhood $V_m$ of $z$ such that
\[
x\in V_m\quad\text{implies}\quad \lambda_x(K_m)<\lambda_z(K_m)+1.
\]
Since $[z]$ is locally closed, by Lemma~4.1(1) in \cite{Clark-anHuef2010-preprint} we can crop $V_1$ if necessary to ensure that $G_z^{V_1}$ is relatively compact. By further cropping each $V_m$ we may assume that $\{V_m\}$ is a decreasing neighbourhood basis of $z$. By our hypothesis, for each $m$ there exists $n_m$ such that 
\begin{equation}\label{lower_bound}
n\ge n_m\quad\text{implies}\quad\lambda_{x_n}(G^{V_m})>k\big(\lambda_z(K_m)+1\big).
\end{equation}
Note that for any $\gamma\in G_{x_n}^{V_m}$ with $n\ge n_m$, we have $r(\gamma)\in V_m$, and so $\lambda_{r(\gamma)}(K_m)<\lambda_z(K_m)+1$. By Haar-system invariance we know that $\lambda_{r(\gamma)}(K_m)=\lambda_{x_n}(K_m\gamma)$, which shows us that
\begin{equation}\label{upper_bound}
\lambda_{x_n}(K_m\gamma)<\lambda_z(K_m)+1.
\end{equation}
If necessary we can increase the elements of $\{n_m\}$ so that it is a strictly increasing sequence.

We now proceed as in the proof of Proposition \ref{AaH_prop4_1_1}. For all $n<n_1$ and $1\le i\le k$ let $\gamma_n^{(i)}=x_n$. For each $n\ge n_1$ there exists a unique number $m(n)$ such that $n_{m(n)}\le n<n_{m(n)+1}$. For the remainder of this proof we will write $m$ instead of $m(n)$ because the specific $n$ will be clear from the context. For each $n\ge n_1$ choose 
\[
\gamma_n^{(1)}\in G_{x_n}^{V_{m}}.
\]
Then by \eqref{lower_bound} and \eqref{upper_bound} we have
\begin{align*}
\lambda_{x_n}(G^{V_m}\backslash K_m\gamma_n^{(1)})&=\lambda_{x_n}(G^{V_m})-\lambda_{x_n}(G^{V_m}\cap K_m\gamma_n^{(1)})\\
&\ge \lambda_{x_n}(G^{V_m})-\lambda_{x_n}(K_m\gamma_n^{(1)})\\
&>k\big(\lambda_z(K_m)+1\big) - \big(\lambda_z(K_m)+1\big)\\
&=(k-1)\big(\lambda_z(K_m)+1\big).
\end{align*}
We can thus choose $\gamma_n^{(2)}\in G_{x_x}^{V_m}\backslash K_m\gamma_n^{(1)}$ for each $n\ge n_1$. This now gives us
\begin{align*}
\lambda_{x_n}&(G^{V_m}\backslash (K_m\gamma_n^{(1)}\cup K_m\gamma_n^{(2)}))\\
&=\lambda_{x_n}(G^{V_m}\backslash K_m\gamma_n^{(1)})-\lambda_{x_n}\big( (G^{V_m}\backslash K_m\gamma_n^{(1)})\cap K_m\gamma_n^{(1)}\big)\\
&\ge \lambda_{x_n}(G^{V_m}\backslash K_m\gamma_n^{(1)})-\lambda_{x_n}(K_m\gamma_n^{(2)})\\
&>(k-1)\big(\lambda_z(K_m)+1\big) - \big(\lambda_z(K_m)+1\big)\\
&=(k-2)\big(\lambda_z(K_m)+1\big).
\end{align*}
Continuing in this manner we can choose 
\[
\gamma_n^{(j)}\in G^{V_m}_{x_n}\backslash \bigg(\bigcup_{i=1}^{j-1}K_m\gamma_n^{(i)}\bigg)
\]
for every $n\ge n_1$ and $j=3,\ldots,k$. The tail of the proof of Proposition \ref{AaH_prop4_1_1} establishes our desired result.
\end{proof}
\end{prop}

\section{Measure ratios and bounds on lower multiplicity}\label{sec_fun}
In this section we show that upper bounds on measure ratios along orbits gives upper bounds on  multiplicities. A subset $S$ of a topological space $X$ is {\em locally closed} if there exist an open set $U$ of $X$ and a closed set $V$ of $X$ such that $S=U\cap V$; this is equivalent to $S$ being open in the closure of $S$ with the subspace topology by, for example, \cite[Lemma~1.25]{Williams2007}.

\begin{lemma}\label{based_on_Ramsay}
Suppose $G$ is a second-countable locally-compact Hausdorff groupoid. Suppose $z\in G^{(0)}$ and $[z]$ is locally closed. Then the restriction of $r$ to $G_z/(G|_{\{z\}})$ is a homeomorphism onto $[z]$. If in addition $G$ is principal, then the restriction of $r$ to $G_z$ is a homeomorphism onto $[z]$.
\begin{proof}
We consider the transitive groupoid $G|_{[z]}$. Since $[z]$ is locally closed, $G|_{[z]}$ is a second-countable locally-compact Hausdorff groupoid. Thus $G|_{[z]}$ is Polish by, for example, \cite[Lemma~6.5]{Williams2007}. Now \cite[Theorem~2.1]{Ramsay1990} applies to give the result.
\end{proof}
\end{lemma}

Theorem~\ref{M2_thm} is based on \cite[Theorem~3.1]{Archbold-anHuef2006}; it is only an intermediary result which will be used to prove a sharper bound in Theorem \ref{M_thm}.
\begin{theorem}\label{M2_thm}
Suppose $G$ is a second-countable locally-compact Hausdorff principal groupoid with Haar system $\lambda$. Let $M\in\RR$ with $M\ge 1$, suppose $z\in G^{(0)}$ such that $[z]$ is locally closed and let $\{x_n\}$ be a sequence in $G^{(0)}$. Suppose there exists an open neighbourhood $V$ of $z$ in $G^{(0)}$ such that $G_z^V$ is relatively compact and 
\[
\lambda_{x_n}(G^V)\le M\lambda_z(G^V)
\]
frequently. Then $\M_\L(\L^z,\{\L^{x_n}\})\le\lfloor M^2\rfloor$.
\begin{proof}[Proofish]
Fix $\epsilon>0$ such that $M^2(1+\epsilon)^2<\lfloor M^2\rfloor +1$. We will build a function $D\in C_c(G)$ such that $\L^z(D^\ast\ast D)$ is a rank-one projection and
\[
\tr\big(\L^{x_n}(D^\ast\ast D)\big)<M^2(1+\epsilon)^2<\lfloor M^2\rfloor +1
\]
frequently. By the generalised lower semi-continuity result of \cite[Theorem~4.3]{Archbold-Spielberg1996} we will have
\begin{align*}
\liminf\,\tr\big(\L^{x_n}(D^\ast\ast D)\big)&\ge \M_\L(\L^z,\{\L^{x_n}\})\, \tr \big(\L^z(D^*\ast D)\big)\\
&=\M_\L(\L^z,\{\L^{x_n}\}),
\end{align*}
and the result will follow.

For the next few paragraphs we will be working with $G_z$ equipped with the subspace topology. Note that $\lambda_z$ can be thought of as a Radon measure on $G_z$ with $\lambda_z(S\cap G_z)=\lambda_z(S)$ for any $\lambda_z$-measurable subset $S$ of $G$. Fix $\delta>0$ such that
\[
\delta<\frac{\epsilon\lambda_z(G^V)}{1+\epsilon}<\lambda_z(G^V).
\]
Since $\lambda_z$ is inner regular on open sets and $G_z^V$ is $G_z$-open, there exists a $G_z$-compact subset $W$ of $G_z^V$ such that
\[
0<\lambda_z(G_z^V)-\delta<\lambda_z(W).
\]
Since $W$ is $G_z$-compact there exists a $G_z$-compact neighbourhood $W_1$ of $W$ that is contained in $G_z^V$ and there exists a continuous function $g:G_z\rightarrow [0,1]$ that is identically one on $W$ and zero off the interior of $W_1$. We have
\[
\lambda_z(G^V)-\delta=\lambda_z(G_z^V)-\delta<\lambda_z(W)\le\int_{G_z} g(t)^2\, d\lambda_z(t)=\|g\|_z^2,
\]
and hence
\begin{equation}\label{AaH3.1}
\frac{\lambda_z(G^V)}{\|g\|_z^2}<1+\frac{\delta}{\|g\|_z^2}<1+\frac{\delta}{\lambda_z(G^V)-\delta}<1+\epsilon.
\end{equation}

By Lemma \ref{based_on_Ramsay} the restriction $\tilde{r}$ of $r$ to $G_z$ is a homeomorphism onto $[z]$. So there exists a continuous function $g_1:\tilde{r}(W_1)\rightarrow [0,1]$ such that $g_1\big(\tilde{r}(\gamma)\big)=g(\gamma)$ for all $\gamma\in W_1$. Thus $\tilde{r}(W_1)$ is $[z]$-compact, which implies that $\tilde{r}(W_1)$ is $G^{(0)}$-compact. Since we know that $G^{(0)}$ is second countable and Hausdorff, Tietze's Extension Theorem can be applied to extend $g_1$ to a continuous map
\[
g_2:G^{(0)}\rightarrow [0,1].
\]
Because $\tilde{r}(W_1)$ is a compact subset of the open set $V$, there exist a compact neighbourhood $P$ of $\tilde{r}(W_1)$ contained in $V$ and a continuous function $h:G^{(0)}\rightarrow [0,1]$ that is identically one on $\tilde{r}(W_1)$ and zero off the interior of $P$. Note that $h$ has compact support that is contained in $P$.

We set $f(x)=h(x)g_2(x)$. Then $f\in C_c(G^{(0)})$ with $0\le f\le 1$ and 
\begin{equation}\label{supp_f_contained_in_V}
\supp f\subset \supp h\subset P\subset V.
\end{equation}
Note that
{\allowdisplaybreaks
\begin{align}
\|f\circ r\|_z^2&=\int_{G_z} f\big(\tilde{r}(\gamma)\big)^2\,d\lambda_z(\gamma)\notag\\
&=\int_{G_z} h\big(\tilde{r}(\gamma)\big)^2g_2\big(\tilde{r}(\gamma)\big)^2\, d\lambda_z(\gamma)\notag\\
&\ge\int_{W_1} h\big(\tilde{r}(\gamma)\big)^2 g(\gamma)^2\, d\lambda_z(\gamma)\notag\\
&=\int_{W_1}g(\gamma)^2\,d\lambda_z(\gamma)\notag\\
&=\|g\|_z^2\label{AaH3.2}
\end{align}}since $\supp\, g\subset W_1$ and $h$ is identically one on $\tilde{r}(W_1)$. We now define $F\in C_c(G^{(0)})$ by
\begin{equation}\label{definition_F}
F(x)=\frac{f(x)}{\|f\circ r\|_z}.
\end{equation}
Then $\|F\circ r\|_z=1$ and
\begin{equation}\label{unmotivated_ref}
F\circ r(\gamma)\ne 0\implies h\big(r(\gamma)\big)\ne 0 \implies r(\gamma)\in V \implies \gamma\in G^V.
\end{equation}

Let $N=\supp\, F$ so that $N=\supp\, f\subset V$ by \eqref{supp_f_contained_in_V} and \eqref{definition_F}. Since $G_z^V$ is relatively compact by our hypothesis, the set $\overline{G_z^N}$ is compact. Let $b\in C_c(G)$ be a function that is identically one on $(\overline{G_z^N})(\overline{G_z^N})^{-1}$ and has range contained in $[0,1]$. We can assume that $b$ is self-adjoint by considering $\frac12(b+b^*)$ if necessary. Define $D\in C_c(G)$ by
\[
D(\gamma):=F\big(r(\gamma)\big)F\big(s(\gamma)\big)b(\gamma).
\]
For $\xi\in L^2(G,\lambda_u)$ and $\gamma\in G$ we have
\begin{align*}
\big(\L^{u}(D)\xi\big)(\gamma)&=\int_GD(\gamma\alpha^{-1})\xi(\alpha)\,d\lambda_u(\alpha)\\
&=\int_G F\big(r(\gamma)\big)F\big(s(\alpha^{-1})\big) b(\gamma\alpha^{-1})\xi(\alpha)\, d\lambda_u(\alpha)\\
&=F\big(r(\gamma)\big) \int_G F\big(r(\alpha)\big) b(\gamma\alpha^{-1})\xi(\alpha)\, d\lambda_u(\alpha).
\end{align*}

In the case where $u=z$, if $\alpha,\gamma\in \supp\, F\circ r\cap s^{-1}(z)$, then $r(\alpha),r(\gamma)\in \supp\, F=N$ and $\gamma,\alpha\in G_z^N$. This implies $b(\gamma\alpha^{-1})=1$, so
\begin{align*}
\big(\L^z(D)\xi\big)(\gamma)&=\int_G F\big(r(\alpha)\big) \xi(\alpha)\, d\lambda_z(\alpha)\\
&=(\xi\,|\,F\circ r)_z F\circ r(\gamma),
\end{align*}
and $L^z(D)$ is a rank-one projection.

By the hypothesis on $V$ there exists a subsequence $\{x_{n_i}\}$ of $\{x_n\}$ such that
\[
\lambda_{x_{n_i}}(G^V)\le M\lambda_z(G^V)
\]
for all $i\ge 1$. If we define $E:=\{\gamma\in G:F\big(r(\gamma)\big)\ne 0\}$ then $E$ is open with
\begin{align}
\lambda_{x_{n_i}}(E)&\le\lambda_{x_{n_i}}(G^V)\text{\quad(using \eqref{unmotivated_ref})}\notag\\
&\le M\lambda_z(G^V)\label{measure_of_E_finite}
\end{align}
and
\begin{equation}\label{AaH3.3}
\int_G\big(F\circ r(\gamma)\big)^2\, d\lambda_{x_{n_i}}(\gamma)\le \frac{\lambda_{x_{n_i}}(E)}{\|f\circ r\|_z^2}
\le \frac{M\lambda_z(G^V)}{\|g\|_z^2}.
\end{equation}
by \eqref{AaH3.2}. Consider the continuous function
\[
T(\alpha,\beta):=F\big(r(\alpha)\big) F\big(r(\beta)\big) b(\alpha\beta^{-1}).
\]
Note that
\begin{align*}
\int_G &T(\alpha,\beta)^2\, d(\lambda_{x_{n_i}}\times\lambda_{x_{n_i}})(\alpha,\beta)\\
&=\int_G F\big(r(\alpha)\big)^2F\big(r(\beta)\big)^2b(\alpha\beta^{-1})^2\, d(\lambda_{x_{n_i}}\times\lambda_{x_{n_i}})(\alpha,\beta)\\
&\le\|F\|_\infty^4\int_G \chi_{E\times E}(\alpha,\beta)\,d(\lambda_{x_{n_i}}\times\lambda_{x_{n_i}})(\alpha,\beta)\\
&=\|F\|_\infty^4\lambda_{x_{n_i}}(E)^2,
\end{align*}
which is finite by \eqref{measure_of_E_finite}. Thus
\[T\in L^2(G\times G,\lambda_{x_{n_i}}\times \lambda_{x_{n_i}}),\]
and since $T$ is conjugate symmetric, \cite[Proposition~3.4.16]{Pedersen1989} implies that $\L^{x_{n_i}}(D)$ is the self-adjoint Hilbert-Schmidt operator on $L^2(G,\lambda_{x_{n_i}})$ with kernel $T$. It follows that $\L^{x_{n_i}}(D^*\ast D)$ is a trace-class operator, and since we equip the Hilbert-Schmidt operators with the trace norm, we have
\[
\tr\,\L^{x_{n_i}}(D^*\ast D)=\|T\|_{L^2(\lambda_{x_{n_i}}\times\lambda_{x_{n_i}})}^2.
\]
Applying Fubini's Theorem to $T$ now gives
{\allowdisplaybreaks\begin{align}
\tr\,\L^{x_{n_i}}&(D^*\ast D)\notag\\
&=\int_G\int_G F\big(r(\alpha)\big)^2 F\big(r(\beta)\big)^2b(\alpha\beta^{-1})^2\, d\lambda_{x_{n_i}}(\alpha)\, d\lambda_{x_{n_i}}(\beta)\notag\\
&\le\bigg(\int_G F\big(r(\alpha)\big)^2\, d\lambda_{x_{n_i}}(\alpha)\bigg)^2\notag\\
&\le\frac{M^2\lambda_z(G^V)^2}{\|g\|_z^4}\text{\quad (using \eqref{AaH3.3})}\notag\\
&<M^2(1+\epsilon)^2\text{\quad (using \eqref{AaH3.1})}\label{Aah3.4}.
\end{align}
}Now
\begin{align*}
\M_\L(\L^z,\{\L^{x_n}\})&\le \underset{\scriptstyle n}{\lim\,\inf}\, \tr\big(\L^{x_n}(D^*\ast D)\big)\\
&\le M^2(1+\epsilon)^2\\
&<\lfloor M^2\rfloor+1,
\end{align*}
and hence $\M_\L(\L^z,\{\L^{x_n}\})\le\lfloor M^2\rfloor$, completing the proof.
\end{proof}
\end{theorem}

The following proposition is an immediate consequence of Theorem \ref{M2_thm} and Proposition \ref{part_of_AaH4.2_generalisation}. This result will be strengthened later in Corollary \ref{AaH_cor5.6}, where we will show that these three items are in fact equivalent.
\begin{prop}\label{AaH_prop_4.2}
Suppose $G$ is a second-countable locally-compact Hausdorff principal groupoid with Haar system $\lambda$. Let $z\in G^{(0)}$ and let $\{x_n\}$ be a sequence in $G^{(0)}$. Assume that $[z]$ is locally closed in $G^{(0)}$. Consider the following properties.
\begin{enumerate}\renewcommand{\labelenumi}{(\arabic{enumi})}
\item\label{AaH_prop_4_2_1} $\M_\L(\L^z,\{\L^{x_n}\})=\infty$.
\item\label{AaH_prop_4_2_2} For every open neighbourhood $V$ of $z$ such that $G_z^V$ is relatively compact, $\lambda_{x_n}(G^V)\rightarrow\infty$ as $n\rightarrow\infty$.
\item\label{AaH_prop_4_2_3} For each $k\ge 1$, the sequence $\{x_n\}$ converges $k$-times in $G^{(0)}/G$ to $z$.
\end{enumerate}
Then {\rm\eqref{AaH_prop_4_2_1}} implies {\rm\eqref{AaH_prop_4_2_2}} and {\rm\eqref{AaH_prop_4_2_2}} implies {\rm\eqref{AaH_prop_4_2_3}}.
\end{prop}

Our next goal is to sharpen the $\lfloor M^2\rfloor$ bound in Theorem \ref{M2_thm}. This strengthened theorem appears later on as Theorem \ref{M_thm}. We will first establish several results to assist in strengthening this bound.

\begin{lemma}\label{lemma_orbits_equal}
Suppose $G$ is a second-countable groupoid and $x,y\in G^{(0)}$. If $\overline{[x]}=\overline{[y]}$ and $[x]$ is locally closed, then $[x]=[y]$.
\begin{proof}
We have $x\in \overline{[y]}$, so there exists $\{\gamma_n\}\subset G$ such that $s(\gamma_n)=y$ and $r(\gamma_n)\rightarrow x$. Since $[x]$ is locally closed, there exists an open subset $U$ of $G$ such that $[x]=U\cap\overline{[x]}$. Then $r(\gamma_n)$ is eventually in $U$, so eventually $r(\gamma_n)\in U\cap\overline{[y]}=U\cap\overline{[x]}=[x]$. Thus there exists $\gamma\in G$ with $s(\gamma)=y$ and $r(\gamma)\in [x]$, as required.
\end{proof}
\end{lemma}

\begin{lemma}\label{the_unbroken_lemma}
Suppose $G$ is a second-countable  groupoid with Haar system $\lambda$. Let $W$ be a compact neighbourhood of $z\in G^{(0)}$ and let $K$ be a compact subset of $G$. Let $\{x_n\}$ be a sequence in $G^{(0)}$ such that $[x_n]\rightarrow [z]$ uniquely in $G^{(0)}/G$. Then for every $\delta>0$ there exists $n_0$ such that, for every $n\ge n_0$ and every $\gamma\in G_{x_n}^W$,
\[
\lambda_{x_n}(K\gamma\cap G^W)<\lambda_z(G^W)+\delta.
\]
\begin{proof}
Suppose not. Then, by passing to a subsequence if necessary, for each $n$ there exists $\gamma_n\in G_{x_n}^W$ such that
\begin{equation}\label{assumption_in_unbroken_lemma}
\lambda_{x_n}(K\gamma_n\cap G^W)\ge \lambda_z(G^W)+\delta.
\end{equation}
Since each $r(\gamma_n)$ is in the compact set $W$, we can pass to a subsequence so that $r(\gamma_n)\rightarrow y$ for some $y\in G^{(0)}$. This implies $[r(\gamma_n)]\rightarrow [y]$, but $[r(\gamma_n)]=[s(\gamma_n)]=[x_n]$ and $[x_n]\rightarrow [z]$ uniquely, so $[y]=[z]$. Choose $\psi\in G$ with $s(\psi)=z$ and $r(\psi)=y$. By Haar-system invariance
\[
\lambda_{x_n}(K\gamma_n\cap G^W)=\lambda_{r(\gamma_n)}(K\cap G^W),
\]
so by applying Lemma \ref{astrids_lim_sup} with the compact space $K\cap G^W$ and $\{r(\gamma_n)\}$ converging to $y$,
\begin{align*}
\underset{\scriptstyle n}{\lim\,\sup}\,\lambda_{x_n}(K\gamma_n\cap G^W)
&=\underset{\scriptstyle n}{\lim\,\sup}\,\lambda_{r(\gamma_n)}(K\cap G^W)\\
&\le\lambda_y(K\cap G^W)\quad\text{(by Lemma \ref{astrids_lim_sup})}\\
&=\lambda_z(K\psi\cap G^W)\quad\text{(Haar-system invariance)}\\
&\le\lambda_z(G^W).
\end{align*}
This contradicts our assertion \eqref{assumption_in_unbroken_lemma}.
\end{proof}
\end{lemma}

The following is a generalisation of \cite[Lemma~3.3]{Archbold-anHuef2006}.
\begin{lemma}\label{AaH_Lemma3.3}
Suppose $G$ is a groupoid with Haar system $\lambda$. Fix $\epsilon>0$, $z\in G^{(0)}$ and let $V$ be an open neighbourhood of $z\in G^{(0)}$ such that $\lambda_z(G^V)<\infty$. Then there exists an open relatively-compact neighbourhood $V_1$ of $z$ such that $\overline{V_1}\subset V$ and
\[
\lambda_z(G^V)-\epsilon<\lambda_z(G^{V_1})\le\lambda_z(G^{\overline{V_1}})\le\lambda_z(G^V)<\lambda_z(G^{V_1})+\epsilon.
\]
\begin{proof}
We use $G_z$ equipped with the subspace topology to find a compact subset $\lambda_z$-estimate of $V$. This estimate is then used to obtain the required open set $V_1$. Since $G_z^V$ is $G_z$-open, by the regularity of $\lambda_z$ there exists a compact subset $W$ of $G_z^V$ such that $\lambda_z(W)>\lambda_z(G_z^V)-\epsilon$. Then $r(W)$ is compact and contained in $V$, so there exists an open relatively-compact neighbourhood $V_1$ of $r(W)$ such that $\overline{V_1}\subset V$. Then
\begin{align*}
\lambda_z(G^V)-\epsilon<\lambda_z(W)&\le \lambda_z(G^{V_1}) \le \lambda_z(G^{\overline{V_1}})\le \lambda_z(G^V)\\&<\lambda_z(W)+\epsilon\le\lambda_z(G^{V_1})+\epsilon,
\end{align*}
as required.
\end{proof}
\end{lemma}

The following lemma is equivalent to the claim in \cite[Proposition~3.6]{Clark2007} that $[x]\mapsto [L^x]$ from $G^{(0)}/G$ to the spectrum of $C^*(G)$ is open.
\begin{lemma}\label{lemma_ind_reps_converge_imply_orbits_converge}
Suppose $G$ is a second-countable locally-compact Hausdorff groupoid with Haar system $\lambda$. If $\{x_n\}$ is a sequence in $G^{(0)}$ with $\L^{x_n}\rightarrow \L^{z}$, then $[x_n]\rightarrow [z]$.
\begin{proof}
We prove the contrapositive. Suppose $[x_n]\nrightarrow [z]$. Then there exists an open neighbourhood $U_0$ of $[z]$ in $G^{(0)}/G$ such that $[x_n]$ is frequently not in $U_0$. Let $q:G^{(0)}\rightarrow G^{(0)}/G$ be the quotient map $x\mapsto [x]$. Then $U_1:=q^{-1}(U_0)$ is an open invariant neighbourhood of $z$ and $x_n\notin U_1$ frequently. Note that $C^*(G|_{U_1})$ is isomorphic to a closed two-sided ideal $I$ of $C^*(G)$ (see \cite[Lemma~2.10]{Muhly-Renault-Williams1996}). 

We now claim that $I\subset\ker\,\L^{x_n}$ whenever $x_n\notin U_1$. Suppose $x_n\notin U_1$ and recall from Remark \ref{measure_induced_epsilon_x} that $\L^{x_n}$ acts on $L^2(G,\lambda_{x_n})$. Fix $f\in C_c(G)$ such that $f(\gamma)=0$ whenever $\gamma\notin G|_{U_1}$ and fix $\xi\in L^2(G,\lambda_{x_n})$. Then by Remark \ref{measure_induced_epsilon_x} we have
\[
\|\L^{x_n}(f)\xi\|_{x_n}^2=\int_G\bigg(\int_G f(\gamma\alpha^{-1})\xi(\alpha)\, d\lambda_{x_n}(\alpha)\bigg)^2\, d\lambda_{x_n}(\gamma).
\]
When evaluating the inner integrand, we have $s(\alpha)=s(\gamma)=x_n$, so $\gamma\alpha^{-1}\in G|_{[x_n]}$. Since $U_1$ is invariant with $x_n\notin U_1$, it follows that $\gamma\alpha^{-1}\notin G|_{U_1}$, and so $f(\gamma\alpha^{-1})=0$. Thus
\[
\|\L^{x_n}(f)\xi\|_{x_n}=0,
\]
and since $\xi$ was fixed arbitrarily, $\L^{x_n}(f)=0$. This implies that $I\subset \ker\,\L^{x_n}$.

We now conclude by observing that since $I\subset\ker\,\L^{x_n}$ frequently, $\L^{x_n}\notin \hat{I}$ frequently. But $\hat{I}$ is an open neighbourhood of $\L^z$, so $\L^{x_n}\nrightarrow \L^z$.
\end{proof}
\end{lemma}

We may now proceed to strengthening the $\lfloor M^2\rfloor$ bound in Theorem \ref{M2_thm}. This theorem is a generalisation of \cite[Theorem~3.5]{Archbold-anHuef2006}.
\begin{theorem}\label{M_thm}
Suppose $G$ is a second-countable locally-compact Hausdorff principal groupoid with Haar system $\lambda$. Let $M\in\RR$ with $M\ge 1$, suppose $z\in G^{(0)}$ such that $[z]$ is locally closed and let $\{x_n\}$ be a sequence in $G^{(0)}$. Suppose there exists an open neighbourhood $V$ of $z$ in $G^{(0)}$ such that $G_z^V$ is relatively compact and 
\[
\lambda_{x_n}(G^V)\le M\lambda_z(G^V)
\]
frequently. Then $\M_\L(\L^z,\{\L^{x_n}\})\le\lfloor M\rfloor$.
\begin{proof}
If $\L^{x_n}$ does not converge to $\L^z$, then $\M_\L(\L^z,\{\L^{x_n}\})=0<\lfloor M \rfloor$. So we assume from now on that $\L^{x_n}\rightarrow \L^z$. Lemma \ref{lemma_ind_reps_converge_imply_orbits_converge} now shows that $[x_n]\rightarrow [z]$.

Next we claim that we may assume, without loss of generality, that $[z]$ is the unique limit of $\{[x_n]\}$ in $G^{(0)}/G$. To see this, note that $\M_\L(\L^z,\{\L^{x_n}\})\le\lfloor M^2\rfloor<\infty$ by Theorem \ref{M2_thm}. Hence, by \cite[Proposition~3.4]{Archbold-anHuef2006}, $\{\L^z\}$ is open in the set of limits of $\{\L^{x_n}\}$. So there exists an open neighbourhood $U_2$ of $\L^z$ in $C^*(G)^\wedge$ such that $\L^z$ is the unique limit of $\{\L^{x_n}\}$ in $U_2$.

By \cite[Proposition~2.5]{Muhly-Williams1990} there is a continuous function $\L:G^{(0)}/G\rightarrow C^*(G)^\wedge$ such that $[x]\mapsto\L^x$ for all $x\in G^{(0)}$. Define $p:G^{(0)}\rightarrow G^{(0)}/G$ by $p(x)=[x]$ for all $x\in G^{(0)}$. Then $p$ is continuous, and
\[
Y:=(\L\circ p)^{-1}(U_2).
\]
is an open $G$-saturated neighbourhood of $z$ in $G^{(0)}$. Note that $x_n\in Y$ eventually.

Now suppose that, for some $y\in Y$, $[x_n]\rightarrow [y]$ in $Y/G$ and hence in $G^{(0)}/G$. Then $\L^{x_n}\rightarrow \L^y$ by \cite[Proposition~2.5]{Muhly-Williams1990}, and $\L^y\in U_2$ since $y\in(\L\circ p)^{-1}(U_2)$. But $\{\L^{x_n}\}$ has the unique limit $\L^z$ in $U_2$, so $\L^z=\L^y$ and hence $\overline{[z]}=\overline{[y]}$. Since $[z]$ is locally closed, Lemma \ref{lemma_orbits_equal} shows that $[z]=[y]$ in $G^{(0)}$ and hence in $Y$.

We know $Y$ is an open saturated subset of $G^{(0)}$, so $C^*(G|_Y)$ is isomorphic to a closed two-sided ideal $J$ of $C^*(G)$. We can apply \cite[Proposition~5.3]{Archbold-Somerset-Spielberg1997} with the $C^*$-subalgebra $J$ to see that $\M_\L(\L^z,\{\L^{x_n}\})$ is the same whether we compute it in the ideal $J$ or in $C^*(G)$. Since $Y$ is $G$-invariant, $G_z^V=G_z^{V\cap Y}$ and eventually $G_{x_n}^V=G_{x_n}^{V\cap Y}$. We may thus consider $G|_Y$ instead of $G$ and therefore assume that $[z]$ is the unique limit of $[x_n]$ in $G^{(0)}/G$ as claimed.

As in \cite{Archbold-anHuef2006}, the idea for the rest of the proof is the same as in Theorem \ref{M2_thm}, although more precise estimates are used. Fix $\epsilon>0$ such that $M(1+\epsilon)^2<\lfloor M\rfloor +1$ and choose $\kappa>0$ such that
\begin{equation}\label{choice_of_kappa}
\kappa<\frac{\epsilon\lambda_z(G^V)}{1+\epsilon}<\lambda_z(G^V).
\end{equation}
By Lemma \ref{AaH_Lemma3.3} there exists an open relatively compact neighbourhood $V_1$ of $z$ such that $\overline{V_1}\subset V$ and
\[
0<\lambda_z(G^V)-\kappa<\lambda_z(G^{V_1})\le\lambda_z(G^{\overline{V_1}})\le\lambda_z(G^V)<\lambda_z(G^{V_1})+\kappa.
\]
Choose a subsequence $\{x_{n_i}\}$ of $\{x_n\}$ such that
\[
\lambda_{x_{n_i}}(G^V)\le M\lambda_z(G^V)
\]
for all $i\ge 1$. Then
\begin{align}
\lambda_{x_{n_i}}(G^{V_1})&\le\lambda_{x_{n_i}}(G^V)\notag\\
&\le M\lambda_z(G^V)\notag\\
&<M\big(\lambda_z(G^{V_1})+\kappa\big)\notag\\
&<M\lambda_z(G^{V_1})+M\epsilon\big(\lambda_z(G^V)-\kappa\big)\quad\text{(by \eqref{choice_of_kappa})}\notag\\
&<M\lambda_z(G^{V_1})+M\epsilon\lambda_z(G^{V_1})\notag\\
&=M(1+\epsilon)\lambda_z(G^{V_1})\label{AaH3.8}
\end{align}
for all $i$. Since
\[
\frac{\lambda_z(G^{V_1})\big(\lambda_z(G^{V_1})+\kappa+1/j\big)}{\big(\lambda_z(G^{V_1})-1/j\big)^2}
\rightarrow
1+\frac{\kappa}{\lambda_z(G^{V_1})}
<1+\epsilon
\]
as $j\rightarrow\infty$, there exists $\delta>0$ such that $\delta<\lambda_z(G^{V_1})$ and
\begin{equation}\label{AaH3.9}
\frac{\lambda_z(G^{V_1})\big(\lambda_z(G^{\overline{V_1}})+\delta\big)}{\big(\lambda_z(G^{V_1})-\delta\big)^2}
<
\frac{\lambda_z(G^{V_1})\big(\lambda_z(G^{V_1})+\kappa+\delta\big)}{\big(\lambda_z(G^{V_1})-\delta\big)^2}<1+\epsilon.
\end{equation}

We will now construct a function $F\in C_c(G^{(0)})$ with support inside $V_1$. Since $\lambda_z$ is inner regular on open sets and $G_z^{V_1}$ is $G_z$-open, there exists a $G_z$-compact subset $W$ of $G_z^{V_1}$ such that
\[
0<\lambda_z(G_z^{V_1})-\delta<\lambda_z(W).
\]
Since $W$ is $G_z$-compact there exists a $G_z$-compact neighbourhood $W_1$ of $W$ that is contained in $G_z^{V_1}$ and there exists a continuous function $g:G_z\rightarrow [0,1]$ that is identically one on $W$ and zero off the interior of $W_1$. We have
\begin{equation}\label{AaH3.10}
\lambda_z(G^{V_1})-\delta<\lambda_z(W)
\le\int_{G_z} g(t)^2\, d\lambda_z(t)
=\|g\|_z^2,
\end{equation}
By Lemma \ref{based_on_Ramsay} the restriction $\tilde{r}$ of $r$ to $G_z$ is a homeomorphism onto $[z]$. So there exists a continuous function $g_1:\tilde{r}(W_1)\rightarrow [0,1]$ such that $g_1\big(\tilde{r}(\gamma)\big)=g(\gamma)$ for all $\gamma\in W_1$. Thus $\tilde{r}(W_1)$ is $[z]$-compact, which implies that $\tilde{r}(W_1)$ is $G^{(0)}$-compact. Since we know that $G^{(0)}$ is second countable and Hausdorff, Tietze's Extension Theorem can be applied to show that $g_1$ can be extended to a continuous map
\[
g_2:G^{(0)}\rightarrow [0,1].
\]
Because $\tilde{r}(W_1)$ is a compact subset of the open set $V_1$, there exist a compact neighbourhood $P$ of $\tilde{r}(W_1)$ contained in $V_1$ and a continuous function $h:G^{(0)}\rightarrow [0,1]$ that is identically one on $\tilde{r}(W_1)$ and zero off the interior of $P$. Note that $h$ has compact support that is contained in $P$.

We set $f(x)=h(x)g_2(x)$. Then $f\in C_c(G^{(0)})$ with $0\le f\le 1$ and 
\begin{equation}\label{supp_f_contained_in_V_2}
\supp f\subset \supp h\subset P\subset V_1.
\end{equation}
Note that
\begin{align}
\|f\circ r\|_z^2&=\int_{G_z} f\big(\tilde{r}(\gamma)\big)^2\,d\lambda_z(\gamma)\notag\\
&=\int_{G_z} h\big(\tilde{r}(\gamma)\big)^2g_2\big(\tilde{r}(\gamma)\big)^2\, d\lambda_z(\gamma)\notag\\
&\ge\int_{W_1} h\big(\tilde{r}(\gamma)\big)^2 g(\gamma)^2\, d\lambda_z(\gamma)\notag\\
&=\int_{W_1}g(\gamma)^2\,d\lambda_z(\gamma)\notag\\
&=\|g\|_z^2\label{AaH3.11}
\end{align}
since $\supp\, g\subset W_1$ and $h$ is identically one on $\tilde{r}(W_1)$. We now define $F\in C_c(G^{(0)})$ by
\begin{equation}\label{definition_F_2}
F(x)=\frac{f(x)}{\|f\circ r\|_z}.
\end{equation}
Then $\|F\circ r\|_z=1$ and
\begin{equation}\label{unmotivated_ref_2}
F\circ r(\gamma)\ne 0\implies h\big(r(\gamma)\big)\ne 0 \implies r(\gamma)\in V_1 \implies \gamma\in G^{V_1}.
\end{equation}

Let $N=\supp\, F$. Suppose $K$ is an open relatively compact symmetric neighbourhood of $(\overline{G_z^N})(\overline{G_z^N})^{-1}$ in $G$ and choose $b\in C_c(G)$ such that $b$ is identically one on $(\overline{G_z^N})(\overline{G_z^N})^{-1}$ and identically zero off $K$. As in Theorem \ref{M2_thm} we may assume that $b$ is self-adjoint by considering $\frac12(b+b^*)$. Define $D\in C_c(G)$ by $D(\gamma):=F\big(r(\gamma)\big)F\big(s(\gamma)\big)b(\gamma)$. By the same argument as in Theorem \ref{M2_thm}, $\L^z(D)$, and hence $\L^z(D^*\ast D)$, is the rank one projection determined by the unit vector $F\circ r\in L^2(G,\lambda_z)$. From \eqref{Aah3.4} we have
\begin{align*}
\tr\big(\L^{x_{n_i}}&(D^*\ast D)\big)\\
&=\int_G F\big(r(\beta)\big)^2\bigg(\int_G F\big(r(\alpha)\big)^2 b(\alpha\beta^{-1})^2\, d\lambda_{x_{n_i}}(\alpha)\bigg)\, d\lambda_{x_{n_i}}(\beta).
\end{align*}
Since $b$ is identically zero off $K$, the inner integrand is zero unless $\alpha\beta^{-1}\in K$. Combining this with \eqref{supp_f_contained_in_V_2} and the fact that $\supp\,\lambda_{x_{n_i}}\subset G_{x_{n_i}}$ enables us to see that this inner integrand is zero unless $\alpha\in G_{x_{n_i}}^{V_1}\cap K\beta$. Thus
\begin{align*}
\tr\big(\L^{x_{n_i}}&(D^*\ast D)\big)\\
&\le\int_{\beta\in G_{x_{n_i}}^{V_1}} F\big(r(\beta)\big)^2\bigg(\int_{\alpha\in G_{x_{n_i}}^{V_1}\cap K\beta} F\big(r(\alpha)\big)^2\, d\lambda_{x_{n_i}}(\alpha)\bigg)\, d\lambda_{x_{n_i}}(\beta).\\
&\le\frac{1}{\|f\circ r\|_z^4}\int_{\beta\in G_{x_{n_i}}^{V_1}} 1\bigg(\int_{\alpha\in G_{x_{n_i}}^{V_1}\cap K\beta} 1\, d\lambda_{x_{n_i}}(\alpha)\bigg)\, d\lambda_{x_{n_i}}(\beta).
\end{align*}

Since $\overline{V_1}$ and $\overline{K}$ are compact, by Lemma \ref{the_unbroken_lemma} there exists $i_0$ such that for every $i\ge i_0$ and any $\beta\in G_{x_{n_i}}^{\overline{V_1}}$,
\[
\lambda_{x_{n_i}}(K\beta\cap G^{\overline{V_1}})<\lambda_z(G^{\overline{V_1}})+\delta.
\]
So, provided $i\ge i_0$,
{\allowdisplaybreaks\begin{align*}
\tr\big(\L^{x_{n_i}}(D^*\ast D)\big)&\le \frac{1}{\|f\circ r\|_z^4}\int_{\beta\in G_{x_{n_i}}^{V_1}}\lambda_{x_{n_i}}(K\beta\cap G_{x_{n_i}}^{V_1})\, d\lambda_{x_{n_i}}(\beta)\\
&\le\frac{1}{\|f\circ r\|_z^4}\int_{\beta\in G_{x_{n_i}}^{V_1}}\big(\lambda_z(G_z^{\overline{V_1}})+\delta\big)\, d\lambda_{x_{n_i}}(\beta)\\
&<\frac{\big(\lambda_z(G^{\overline{V_1}})+\delta\big)\lambda_{x_{n_i}}(G^{V_1})}{\|f\circ r\|_z^4}\\
&<\frac{M(1+\epsilon)\big(\lambda_z(G^{\overline{V_1}})+\delta\big)\lambda_z(G^{V_1})}{\|g\|_z^4}\quad\text{(by \eqref{AaH3.8} and \eqref{AaH3.11})}\\
&<\frac{M(1+\epsilon)\big(\lambda_z(G^{\overline{V_1}})+\delta\big)\lambda_z(G^{V_1})}{(\lambda_z(G^{V_1})-\delta)^2}\quad\text{(by \eqref{AaH3.10})}\\
&<M(1+\epsilon)^2\quad\text{(by \eqref{AaH3.9})}.
\end{align*}}

We can now make our conclusion as in \cite[Theorem~3.5]{Archbold-anHuef2006}: by generalised lower semi-continuity \cite[Theorem~4.3]{Archbold-Spielberg1996},
\begin{align*}
\underset{\scriptstyle n}{\lim\,\inf}\,\tr\big(\L^{x_n}(D^\ast\ast D)\big)&\ge \M_\L(\L^z,\{\L^{x_n}\})\, \tr \big(\L^z(D^*\ast D)\big)\\
&=\M_\L(\L^z,\{\L^{x_n}\}).
\end{align*}
We now have
\begin{align*}
\M_\L(\L^z,\{\L^{x_n}\})&\le \underset{\scriptstyle n}{\lim\,\inf}\, \tr\big(\L^{x_n}(D^*\ast D)\big)\\
&\le M(1+\epsilon)^2\\
&<\lfloor M\rfloor+1,
\end{align*}
and so $\M_\L(\L^z,\{\L^{x_n}\})\le\lfloor M\rfloor$, as required.
\end{proof}
\end{theorem}

\section{Lower multiplicity and $k$-times convergence II}\label{sec_lower_multiplicity_2}
We proved in Proposition \ref{AaH_thm_1.1_1_implies_2} that if a sequence converges $k$-times in the orbit space of a principal groupoid, then the lower multiplicity of the associated sequence of representations is at least $k$. In this section we will prove the converse.

Our next lemma generalises \cite[Lemma~5.1]{Archbold-anHuef2006}; with the exception of notation changes, the proof is the same as the proof in \cite{Archbold-anHuef2006}.
\begin{lemma}\label{AaH_lemma5.1}
Suppose $G$ is a second-countable locally-compact Hausdorff principal groupoid. Let $k\in\PP$, $z\in G^{(0)}$, and $\{x_n\}$ be a sequence in $G^{(0)}$. Assume that $[z]$ is locally closed in $G^{(0)}$ and that there exists $R>k-1$ such that for every open neighbourhood $U$ of $z$ with $G_z^U$ relatively compact we have
\[
\underset{\scriptstyle n}{\lim\,\inf}\,\lambda_{x_n}(G^U)\ge R\lambda_z(G^U).
\]
Given an open neighbourhood $V$ of $z$ such that $G_z^V$ is relatively compact, there exists a compact neighborhood $N$ of $z$ with $N\subset V$ such that
\[
\underset{\scriptstyle n}{\lim\,\inf}\,\lambda_{x_n}(G^N)>(k-1)\lambda_z(G^N).
\]
\begin{proof}
Apply Lemma \ref{AaH_Lemma3.3} to $V$ with $0<\epsilon<\frac{R-k+1}R \lambda_z(G^V)$ to get an open relatively-compact neighbourhood $V_1$ of $z$ with $\overline{V_1}\subset V$ and
\[
\lambda_z(G^V)-\epsilon<\lambda_z(G^{V_1})\le\lambda_z(G^{\overline{V_1}})\le\lambda_z(G^V)<\lambda_z(G^{V_1})+\epsilon.
\]
Since $G_z^{V_1}$ is relatively compact we have
\begin{align*}
\underset{\scriptstyle n}{\lim\,\inf}\,\lambda_{x_n}(G^{\overline{V_1}})&\ge\underset{\scriptstyle n}{\lim\,\inf}\,\lambda_{x_n}(G^{V_1})\\
&\ge R\lambda_z(G^{V_1})&\text{(by hypothesis)}\\
&>R\big(\lambda_z(G^V)-\epsilon\big)\\
&>(k-1)\lambda_z(G^V)&\text{(by our choice of }\epsilon\text{)}\\
&\ge (k-1)\lambda_z(G^{\overline{V_1}}).
\end{align*}
So we may take $N=\overline{V_1}$.
\end{proof}
\end{lemma}
\begin{remark}\label{AaH_lemma5.1_variant}
The preceding Lemma also holds when $\lim\,\inf$ is replaced by $\lim\,\sup$.  No modification of the proof is needed beyond replacing the two occurrences of $\lim\,\inf$ with $\lim\,\sup$.
\end{remark}

We may now proceed to our main theorem.
\begin{theorem}\label{circle_thm}
Suppose $G$ is a second-countable locally-compact Hausdorff principal groupoid that admits a Haar system $\lambda$. Let $k$ be a positive integer, let $z\in G^{(0)}$ and let $\{x_n\}$ be a sequence in $G^{(0)}$. Assume that $[z]$ is locally closed in $G^{(0)}$. Then the following are equivalent:
\begin{enumerate}\renewcommand{\labelenumi}{(\arabic{enumi})}
\item\label{circle_thm_1} the sequence $\{x_n\}$ converges $k$-times in $G^{(0)}/G$ to $z$;
\item\label{circle_thm_2} $\M_\L(\L^z,\{\L^{x_n}\})\ge k$;
\item\label{circle_thm_3} for every open neighbourhood $V$ of $z$ in $G^{(0)}$ such that $G_z^V$ is relatively compact we have
\[
\underset{\scriptstyle n}{\lim\,\inf}\,\lambda_{x_n}(G^V)\ge k\lambda_z(G^V);
\]
\item\label{circle_thm_4} there exists a real number $R>k-1$ such that for every open neighbourhood $V$ of $z$ in $G^{(0)}$ with $G_z^V$ relatively compact we have
\[
\underset{\scriptstyle n}{\lim\,\inf}\,\lambda_{x_n}(G^V)\ge R\lambda_z(G^V);\quad\text{and}
\]
\item\label{circle_thm_5} there exists a basic decreasing sequence of compact neighbourhoods $\{W_m\}$ of $z$ in $G^{(0)}$ such that, for each $m\ge 1$,
\[
\underset{\scriptstyle n}{\lim\,\inf}\,\lambda_{x_n}(G^{W_m})>(k-1)\lambda_z(G^{W_m}).
\]
\end{enumerate}
\begin{proof}


We know that \eqref{circle_thm_1} implies \eqref{circle_thm_2} by Proposition \ref{AaH_thm_1.1_1_implies_2}.

Suppose \eqref{circle_thm_2}. If $\M_\L(\L^z,\{\L^{x_n}\})\ge k$, then $\M_\L(\L^z,\{\L^{x_n}\})>\lfloor k-\epsilon\rfloor$ for all $\epsilon>0$. By Theorem \ref{M_thm}, for every $G^{(0)}$-open neighborhood $V$ of $z$ such that $G_z^V$ is relatively compact,
\[
\lambda_{x_n}(G^V)>(k-\epsilon)\lambda_z(G^V)
\]
eventually, and hence \eqref{circle_thm_3} holds.

It is immediately true that \eqref{circle_thm_3} implies \eqref{circle_thm_4}.

Suppose \eqref{circle_thm_4}. We will construct the sequence $\{W_m\}$ of compact neighbourhoods inductively. Let $\{V_j\}$ be a basic decreasing sequence of open neighborhoods of $z$ such that $G_z^{V_1}$ is relatively compact (such neighborhoods exist by \cite[Lemma~4.1(1)]{Clark-anHuef2010-preprint}). By Lemma \ref{AaH_lemma5.1} there exists a compact neighbourhood $W_1$ of $z$ such that $W_1\subset V_1$ and
\[
\lambda_{x_n}(G^{W_1})>(k-1)\lambda_z(G^{W_1}).
\]
Now assume there are compact neighbourhoods $W_1,W_2,\ldots,W_m$ of $z$ with $W_1\supset W_2\supset\cdots\supset W_m$ such that
\begin{equation}\label{circle_thm_4_imply_5_eqn}
W_i\subset V_i\quad\text{and}\quad\lambda_{x_n}(G^{W_i})>(k-1)\lambda_z(G^{W_i})
\end{equation}
for all $1\le i\le m$. Apply Lemma \ref{AaH_lemma5.1} to $(\mathrm{Int}\, m)\cap V_{m+1}$ to obtain a compact neighbourhood $W_{m+1}$ of $z$ such that $W_{m+1}\subset(\mathrm{Int}\, W_m)\cap V_{m+1}$ and \eqref{circle_thm_4_imply_5_eqn} holds for $i=m+1$, establishing \eqref{circle_thm_5}.

Suppose \eqref{circle_thm_5}. We begin by showing that $[x_n]\rightarrow [z]$ in $G^{(0)}/G$. Let $q:G^{(0)}\rightarrow G^{(0)}/G$ be the quotient map. Let $U$ be a neighbourhood of $[z]$ in $G^{(0)}/G$ and $V=q^{-1}(U)$. There exists $m$ such that $W_m\subset V$. Since $\lim\,\inf_n\,\lambda_{x_n}(G^{W_m})>0$ there exists $n_0$ such that $G_{x_n}^{W_m}\ne\emptyset$ for all $n\ge n_0$. Thus, for $n\ge n_0$,
\[
[x_n]=q(x_n)\in q(W_m)\subset q(V)=U.
\]
Thus $[x_n]$ is eventually in every neighbourhood of $[z]$ in $G^{(0)}/G$.

Now suppose that $\M_\L(\L^z,\{\L^{x_n}\})<\infty$. Then, as in the proof of Theorem \ref{M_thm}, we may localise to an open invariant neighbourhood $Y$ of $z$ such that $[z]$ is the unique limit in $Y/G$ of $[x_n]$. Eventually $W_m\subset Y$, and so the sequence $\{x_n\}$ converges $k$-times in $Y/(G|_Y)=Y/G$ to $z$ by Proposition \ref{AaH_prop4_1_1} applied to the groupoid $G|_Y$. This implies that the sequence $\{x_n\}$ converges $k$-times in $G^{(0)}/G$.

Finally, if $\M_\L(\L^z,\{\L^{x_n}\})=\infty$, then $\{x_n\}$ converges $k$-times in $G^{(0)}/G$ to $z$ by Proposition \ref{AaH_prop_4.2}, establishing \eqref{2nd_circle_thm_1} and completing the proof.
\end{proof}
\end{theorem}

\begin{cor}
Suppose that $G$ is a second-countable locally-compact Hausdorff principal groupoid such that all the orbits are locally closed. Let $k\in\PP$ and let $z\in G^{(0)}$ such that $[z]$ is not open in $G^{(0)}$. Then the following are equivalent:
\begin{enumerate}\renewcommand{\labelenumi}{(\arabic{enumi})}
\item\label{AaH_cor5.5_1} whenever $\braces{x_n}$ is a sequence in $G^{(0)}$ which converges to $z$ with $[x_n]\ne [z]$ eventually, then $\braces{x_n}$ is $k$-times convergent in $G^{(0)}/G$ to $z$;
\item\label{AaH_cor5.5_2} $\M_\L(\L^z)\ge k$.
\end{enumerate}
\begin{proof}
Assume \eqref{AaH_cor5.5_1}. We must first establish that $\braces{\L^z}$ is not open in $C^*(G)^\wedge$. If this is not the case, then $\braces{\L^z}$ is open and we can apply \cite[Proposition~3.6]{Clark2007} to see that $\braces{[z]}$ is open in $G^{(0)}/G$, and so $[z]$ is open in $G^{(0)}$, contradicting our assumption. Since $\braces{\L^z}$ is not open in $C^*(G)^\wedge$, we can apply \cite[Lemma~A.2]{Archbold-anHuef2006} to see that there exists a sequence $\braces{\pi_i}$ of irreducible representations of $C^*(G)$ such that each $\pi_i$ is not unitarily equivalent to $\L^z$, $\pi_i\rightarrow \L^z$ in $C^*(G)^\wedge$, and
\begin{equation}\label{AaH_5.3}
\M_\L(\L^z)=\M_\L(\L^z,\braces{\pi_i})=\M_\U(\L^z,\braces{\pi_i}).
\end{equation}
Since the orbits are locally closed, the map $G^{(0)}/G\rightarrow C^*(G)^\wedge$ such that $[x]\mapsto \L^x$ is a homeomorphism by \cite[Proposition~5.1]{Clark2007}\footnote{Proposition~5.1 in \cite{Clark2007} states that if a principal groupoid has locally closed orbits, then the map from $G^{(0)}/G$ to $C^*(G)^\wedge$ where $[x]\mapsto \L^x$ is a `homeomorphism from $G^{(0)}/G$ into $C^*(G)^\wedge$'. The proof explicitely shows that this map is a surjection.}. It follows that the mapping $G^{(0)}\rightarrow C^*(G)^\wedge$ such that $x\mapsto \L^x$ is an open surjection, so by \cite[Proposition~1.15]{Williams2007} there is a sequence $\braces{x_n}$ in $G^{(0)}$ such that $x_n \rightarrow z$ and $\braces{\L^{x_n}}$ is unitarily equivalent to a subsequence of $\braces{\pi_i}$.
By \eqref{AaH_5.3}, 
\[
\M_\L(\L^z)=\M_\U(\L^z,\braces{\pi_i})\ge \M_\U(\L^z,\braces{\L^{x_n}})\ge \M_\L(\L^z,\braces{\L^{x_n}}).
\]
We know by \eqref{AaH_cor5.5_1} that $\braces{x_n}$ converges $k$-times to $z$ in $G^{(0)}/G$, so it follows from Theorem \ref{circle_thm} that $\M_\L(\L^z)\ge\M_\L(\L^z,\braces{\L^{x_n}})\ge k$.

Assume \eqref{AaH_cor5.5_2}. If $\braces{x_n}$ is a sequence in $G^{(0)}$ which converges to $z$ such that $[x_n]\ne [z]$ eventually, then
\[
\M_\L(\L^z,\braces{\L^{x_n}})\ge\M_\L(\L^z)\ge k.
\]
By Theorem \ref{circle_thm}, $\braces{x_n}$ is $k$-times convergent to $z$ in $G^{(0)}/G$.
\end{proof}
\end{cor}

The next corollary improves Proposition \ref{AaH_prop_4.2} and is an immediate consequence of Proposition \ref{AaH_prop_4.2} and Theorem \ref{circle_thm}.
\begin{cor}\label{AaH_cor5.6}
Suppose that $G$ is a second-countable locally-compact Hausdorff principal groupoid with Haar system $\lambda$. Let $z\in G^{(0)}$ and let $\braces{x_n}$ be a sequence in $G^{(0)}$. Assume that $[z]$ is locally closed. Then the following are equivalent:
\begin{enumerate}\renewcommand{\labelenumi}{(\arabic{enumi})}
\item\label{AaH_cor5.6_1} $\M_\L(\L^z,\{\L^{x_n}\})=\infty$.
\item\label{AaH_cor5.6_2} For every open neighbourhood $V$ of $z$ such that $G_z^V$ is relatively compact, $\lambda_{x_n}(G^V)\rightarrow\infty$ as $n\rightarrow\infty$.
\item\label{AaH_cor5.6_3} For each $k\ge 1$, the sequence $\{x_n\}$ converges $k$-times in $G^{(0)}/G$ to $z$.
\end{enumerate}
\end{cor}

\section{Upper multiplicity and $k$-times convergence}\label{sec_upper_multiplicity}
The results in this section  are corollaries of Theorems~\ref{M_thm} and~\ref{circle_thm}: they relate $k$-times convergence, measure ratios and upper multiplicity numbers, generalising all the upper-multiplicity results of \cite{Archbold-anHuef2006}. We begin with the upper-multiplicity analogue of Theorem~\ref{M_thm}.
\begin{theorem}\label{thm_AaH3.6}
Suppose that $G$ is a second-countable locally-compact Hausdorff principal groupoid with Haar system $\lambda$. Let $M\in \RR$ with $M\ge 1$, let $z\in G^{(0)}$ and let $\braces{x_n}$ be a sequence in $G^{(0)}$. Assume that $[z]$ is locally closed. Suppose that there exists an open neighbourhood $V$ of $z$ in $G^{(0)}$ such that $G_z^V$ is relatively compact and
\[
\lambda_{x_n}(G^V)\le M\lambda_z(G^V)<\infty
\]
eventually. Then $\M_\U(\L^z,\braces{\L^{x_n}})\le \lfloor M\rfloor$.
\begin{proof}
Since $G$ is second countable, $C^*(G)$ is separable. By \cite[Lemma~A.1]{Archbold-anHuef2006} there exists a sequence $\braces{\L^{x_{n_i}}}$ such that
\[
\M_\U(\L^z,\braces{\L^{x_n}})=\M_\U(\L^z,\braces{\L^{x_{n_i}}})=\M_\L(\L^z,\braces{\L^{x_{n_i}}}).
\]
By Theorem \ref{M_thm}, $\M_\L(\L^z,\braces{\L^{x_{n_i}}})\le \lfloor M\rfloor$, so $\M_\U(\L^z,\braces{\L^{x_n}})\le\lfloor M\rfloor$.
\end{proof}
\end{theorem}

\begin{cor}\label{AaH_cor3.7}
Suppose that $G$ is a second-countable locally-compact Hausdorff principal groupoid with Haar system $\lambda$ such that all the orbits are locally closed. Let $M\in \RR$ with $M\ge 1$ and let $z\in G^{(0)}$. If for every sequence $\braces{x_n}$ in $G^{(0)}$ which converges to $z$ there exists an open neighbourhood $V$ of $z$ in $G^{(0)}$ such that $G_z^V$ is relatively compact and 
\[
\lambda_{x_n}(G^V)\le M\lambda_z(G^V)<\infty
\]
frequently, then $\M_\U(\L^z)\le\lfloor M\rfloor$.
\begin{proof}
Since $G$ is second countable, $C^*(G)$ is separable, and so we can apply \cite[Lemma~1.2]{Archbold-Kaniuth1999} to see that there exists a sequence $\braces{\pi_n}$ in $C^*(G)^\wedge$ that converges to $\L^z$ such that
\[
\M_\L(\L^z,\braces{\pi_n})=\M_\U(\L^z,\braces{\pi_n})=\M_\U(\L^z).
\]
Since the orbits are locally closed, the map $G^{(0)}/G\rightarrow C^*(G)^\wedge$ such that $[x]\mapsto \L^x$ is a homeomorphism by \cite[Proposition~5.1]{Clark2007}. In particular, the mapping $G^{(0)}\rightarrow C^*(G)^\wedge$ such that $x\mapsto \L^x$ is an open surjection, so by \cite[Proposition~1.15]{Williams2007} there exists a sequence $\braces{x_i}$ in $G^{(0)}$ converging to $z$ such that $\braces{[\L^{x_i}]}$ is a subsequence of $\braces{[\pi_n]}$. By Theorem \ref{M_thm}, $\M_\L(\L^z,\braces{\L^{x_n}})\le\lfloor M\rfloor$. Since 
\begin{align*}
\M_\U(\L^z)&=\M_\L(\L^z,\braces{\pi_n})\le\M_\L(\L^z,\braces{\L^{x_i}})\le\M_\U(\L^z,\braces{\L^{x_i}})\\
&\le\M_\U(\L^z,\braces{\pi_n})=\M_\U(\L^z),
\end{align*}
we obtain $\M_\U(\L^z)\le\lfloor M\rfloor$, as required.
\end{proof}
\end{cor}

In Proposition~\ref{AaH_prop4_1_1} we generalised the first part of \cite[Proposition~4.1]{Archbold-anHuef2006}. We will now generalise the second part. The argument we use is similar to that used in Proposition \ref{AaH_prop4_1_1}.
\begin{prop}\label{AaH_prop4_1_2}
Let $G$ be a second-countable locally-compact Hausdorff principal groupoid with Haar system $\lambda$. Let $k\in\PP$ and $z\in G^{(0)}$ with $[z]$ locally closed in $G^{(0)}$. Assume that $\{x_n\}$ is a sequence in $G^{(0)}$ such that $[x_n]\rightarrow [z]$ uniquely in $G^{(0)}/G$. Suppose $\{W_m\}$ is a basic decreasing sequence of compact neighbourhoods of $z$ such that each $m$ satisfies
\[
\underset{\scriptstyle n}{\lim\,\sup}\,\lambda_{x_n}(G^{W_m})>(k-1)\lambda_z(G^{W_m}).
\]
Then there exists a subsequence of $\{x_n\}$ which converges $k$-times in $G^{(0)}/G$ to $z$. 

\begin{proof}
Let $\{K_m\}$ be an increasing sequence of compact subsets of $G$ such that $G=\bigcup_{m\ge 1}\mathrm{Int}\,K_m$. By the regularity of $\lambda_z$, for each $m\ge 1$ there exist $\delta_m>0$ and an open neighbourhood $U_m$ of $G_z^{W_m}$ such that
\begin{equation}\label{AaH4.4}
\underset{\scriptstyle n}{\lim\,\sup}\,\lambda_{x_n}(G^{W_m})>(k-1)\lambda_z(U_m) + \delta_m.
\end{equation}
We will construct, by induction, a strictly increasing sequence of positive integers $\{i_m\}$ such that, for all $m$,
\begin{align}
&\lambda_{x_{i_m}}(K_m\alpha\cap G^{W_m})<\lambda_z(U_m)+\delta_m/k\quad\text{for all }\alpha\in G_{x_{i_m}}^{W_m},\quad\text{and}\label{AaH4.5}\\
&\lambda_{x_{i_m}}(G^{W_m})>(k-1)\lambda_z(U_m)+\delta_m.\label{AaH4.6}
\end{align}

By Lemma \ref{the_unbroken_lemma} with $\delta=\lambda_z(U_1)-\lambda_z(G^{W_1})+\delta_1/k$, there exists $n_1$ such that $n\ge n_1$ implies

\[
\lambda_{x_n}(K_1\alpha\cap G^{W_1})<\lambda_z(U_1)+\delta_1/k\quad\text{for all }\alpha\in G_{x_n}^{W_m}.
\]
By considering \eqref{AaH4.4} with $m=1$ we can choose $i_1\ge n_1$ such that
\[
\lambda_{x_{i_1}}(G^{W_1})>(k-1)\lambda_z(U_1)+\delta_1.
\]
Assuming that $i_1<i_2<\cdots<i_{m-1}$ have been chosen, we can apply Lemma \ref{the_unbroken_lemma} with $\delta=\lambda_z(U_m)-\lambda_z(G^{W_m})+\delta_m/k$ to obtain $n_m>i_{m-1}$ such that
\[
n\ge n_m\quad\text{implies}\quad\lambda_{x_n}(K_m\alpha\cap G^{W_m})<\lambda_z(U_m)+\delta_m/k\quad\text{for all }\alpha\in G_{x_n}^{W_m},
\]
and then by \eqref{AaH4.4} we can choose $i_m\ge n_m$ such that
\[
\lambda_{x_{i_m}}(G^{W_m})>(k-1)\lambda_z(U_m)+\delta_m.
\]

For each $m\in\PP$ choose $\gamma_{i_m}^{(1)}\in G_{x_{i_m}}^{W_m}$ (which is non-empty by \eqref{AaH4.6}). By \eqref{AaH4.5} and \eqref{AaH4.6} we have
\begin{align*}
\lambda_{x_{i_m}}(G^{W_m}\backslash K_m\gamma_{i_m}^{(1)})&=\lambda_{x_{i_m}}(G^{W_m})-\lambda_{x_{i_m}}(G^{W_m}\cap K_m\gamma_{i_m}^{(1)})\\
&>(k-1)\lambda_z(U_m)+\delta_m-\big(\lambda_z(U_m)+\delta_m/k\big)\\
&=(k-2)\lambda_z(U_m)+\frac{k-1}k\delta_m.
\end{align*}
So we can choose $\gamma_{i_m}^{(2)}\in G_{x_{i_m}}^{W_m}\backslash K_m\gamma_{i_m}^{(1)}$. This implies, as in the proof of Proposition \ref{AaH_prop4_1_1}, that
\[
\lambda_{x_{i_m}}\big(G^{W_m}\backslash (K_m\gamma_{i_m}^{(1)}\cup K_m\gamma_{i_m}^{(2)})\big)>(k-3)\lambda_z(U_m)+\frac{(k-2)}k\delta_m,
\]
enabling us to choose $\gamma_{i_m}^{(3)}\in G_{x_{i_m}}^{W_m}\backslash (K_m\gamma_{i_m}^{(1)}\cap K_m\gamma_{i_m}^{(2)})$. Continuing in this way for  $j=3,\ldots,k$, for each $i_m$ we choose
\begin{equation}\label{eqn_choosing_gammas-2}
\gamma_{i_m}^{(j)}\in G_{x_{i_m}}^{W_m}\backslash\bigg(\bigcup_{l=1}^{j-1}K_m\gamma_{i_m}^{(l)}\bigg).
\end{equation}
Note that $\gamma_{i_m}^{(j)}\notin K_m\gamma_{i_m}^{(l)}$ for $1\le l<j\le k$.

We claim that $r(\gamma_{i_m}^{(l)})\rightarrow z$ as $m\rightarrow\infty$ for $1\le l\le k$. To see this, fix $l$ and let $V$ be an open neighbourhood of $z$. Since $\{W_m\}$ is a decreasing neighbourhood basis for $z$ there exists $m_0$ such that $m\ge m_0$ implies $W_m\subset V$, and so $r(\gamma_{i_m}^{(l)})\in W_m\subset V$.

Finally we claim that $\gamma_{i_m}^{(j)}(\gamma_{i_m}^{(l)})^{-1}\rightarrow\infty$ as $m\rightarrow\infty$ for $1\le l<j\le k$. Fix $l<j$ and let $K$ be a compact subset of $G$. There exists $m_0$ such that $K\subset K_m$ for all $m\ge m_0$. By \eqref{eqn_choosing_gammas-2} we know
\begin{align*}
\gamma_{i_m}^{(j)}&\in G_{x_{i_m}}^{W_m}\backslash (K_m\gamma_{i_m}^{(l)})\\
&=\big(G_{x_{i_m}}^{W_m}(\gamma_{i_m}^{(l)})^{-1}\gamma_{i_m}^{(l)}\big)\backslash (K_m\gamma_{i_m}^{(l)})\\
&=\big((G_{x_{i_m}}^{W_m}(\gamma_{i_m}^{(l)})^{-1})\backslash K_m\big)\gamma_{i_m}^{(l)}.
\end{align*}
So provided $m\ge m_0$, $\gamma_{i_m}^{(j)}(\gamma_{i_m}^{(l)})^{-1}\in \big(G_{x_{i_m}}^{W_m}(\gamma_{i_m}^{(l)})^{-1}\big)\backslash K_m\subset G\backslash K_m\subset G\backslash K$, enabling us to conclude that $\{x_{i_m}\}$ converges $k$-times in $G^{(0)}/G$ to $z$.
\end{proof}
\end{prop}

\begin{theorem}\label{2nd_circle_thm}
Suppose that $G$ is a second-countable locally-compact Hausdorff principal groupoid with Haar system $\lambda$. Let $k\in\PP$, let $z\in G^{(0)}$, and let $\braces{x_n}$ be a sequence in $G^{(0)}$ such that $[x_n]$ converges to $[z]$ in $G^{(0)}/G$. Assume that $[z]$ is locally closed. Then the following are equivalent:
\begin{enumerate}\renewcommand{\labelenumi}{(\arabic{enumi})}
\item\label{2nd_circle_thm_1} there exists a subsequence $\braces{x_{n_i}}$ of $\braces{x_n}$ which converges $k$-times in $G^{(0)}/G$ to $z$;
\item\label{2nd_circle_thm_2} $\M_\U(\L^z,\braces{\L^{x_n}})\ge k$;
\item\label{2nd_circle_thm_3} for every open neighbourhood $V$ of $z$ such that $G_z^V$ is relatively compact we have
\[
\underset{\scriptstyle n}{\lim\,\sup}\,\lambda_{x_n}(G^V)\ge k\lambda_z(G^V);
\]
\item\label{2nd_circle_thm_4}
there exists a real number $R>k-1$ such that for every open neighbourhood $V$ of $z$ in $G^{(0)}$ with $G_z^V$ relatively compact we have
\[
\underset{\scriptstyle n}{\lim\,\sup}\,\lambda_{x_n}(G^V)\ge R\lambda_z(G^V);\quad\text{and}
\]
\item\label{2nd_circle_thm_5} there exists a basic decreasing sequence of compact neighbourhoods $\{W_m\}$ of $z$ in $G^{(0)}$ such that, for each $m\ge 1$,
\[
\underset{\scriptstyle n}{\lim\,\sup}\,\lambda_{x_n}(G^{W_m})>(k-1)\lambda_z(G^{W_m}).
\]
\end{enumerate}
\begin{proof}
If \eqref{2nd_circle_thm_1} holds then $\M_\L(\L^z,\braces{\L^{x_{n_i}}})\ge k$ by Theorem \ref{circle_thm}, and so
\[
\M_\U(\L^z,\braces{\L^{x_n}}\ge\M_\U(\L^z,\braces{\L^{x_{n_i}}})\ge\M_\L(\L^z,\braces{\L^{x_{n_i}}})\ge k.
\]

If \eqref{2nd_circle_thm_2} holds then by \cite[Lemma~A.1]{Archbold-anHuef2006} there is a subsequence $\braces{x_{n_r}}$ such that $\M_\L(\L^z,\braces{\L^{x_{n_r}}})=\M_\U(\L^z,\braces{\L^{x_n}})$ so that $\M_\L(\L^z,\braces{\L^{x_{n_r}}})\ge k$. Let $V$ be any open neighbourhood of $z$ in $G^{\bkt{0}}$ such that $G_z^V$ is relatively compact. Then
\[
\underset{\scriptstyle n}{\lim\,\sup}\,\lambda_{x_n}\bkt{G^V}\ge\underset{\scriptstyle r}{\lim\,\sup}\,\lambda_{x_{n_r}}\bkt{G^V}\ge \underset{\scriptstyle r}{\lim\,\inf}\,\lambda_{x_{n_r}}\bkt{G^V}\ge k\lambda_z\bkt{G^V},
\]
using Theorem \ref{circle_thm} for the last step.

That \eqref{2nd_circle_thm_3} implies \eqref{2nd_circle_thm_4} is immediate.

That \eqref{2nd_circle_thm_4} implies \eqref{2nd_circle_thm_5} follows by making references to Remark \ref{AaH_lemma5.1_variant} rather than Lemma \ref{AaH_lemma5.1} in the \eqref{circle_thm_4} implies \eqref{circle_thm_5} component of the proof of Theorem \ref{circle_thm}.

Assume \eqref{2nd_circle_thm_5}. First suppose that $\M_\L(\L^z,\braces{\L^{x_n}})<\infty$. Since $[x_n]\rightarrow[z]$, we can use an argument found at the beginning of the proof of Theorem \ref{M_thm} to obtain an open $G$-invariant neighborhood $Y$ of $z$ in $G^{(0)}$ so that if we define $H:=G|_Y$, there exists a subsequence $\braces{x_{n_i}}$ of $\braces{x_n}$ such that $[x_{n_i}]\rightarrow [z]$ uniquely in $H^{(0)}/H$. Proposition \ref{AaH_prop4_1_2} now shows us that there exists a subsequence $\braces{x_{n_{i_j}}}$ of $\braces{x_{n_i}}$ that converges $k$-times in $H^{(0)}/H$ to $z$. It follows that $\braces{x_{n_{i_j}}}$ converges $k$-times in $G^{(0)}/G$ to $z$.

When $\M_\L(\L^z,\braces{\L^{x_n}})=\infty$, $\braces{x_n}$ converges $k$-times in $G^{(0)}/G$ to $z$  by Corollary \ref{AaH_cor5.6}, establishing \eqref{2nd_circle_thm_1}.
\end{proof}
\end{theorem}

\begin{cor}\label{AaH_cor5.4}
Suppose that $G$ is a second-countable locally-compact Hausdorff principal groupoid such that all the orbits are locally closed. Let $k\in\PP$ and let $z\in G^{(0)}$. Then the following are equivalent:
\begin{enumerate}\renewcommand{\labelenumi}{(\arabic{enumi})}
\item\label{AaH_cor5.4_1} there exists a sequence $\braces{x_n}$ in $G^{(0)}$ which is $k$-times convergent in $G^{(0)}/G$ to $z$;
\item\label{AaH_cor5.4_2} $\M_\U(\L^z)\ge k$.
\end{enumerate}
\begin{proof}
Assume \eqref{AaH_cor5.4_1}. By the definitions of upper and lower multiplicity, 
\[
\M_\U(\L^z)\ge\M_\U(\L^z,\braces{\L^{x_n}})\ge\M_\L(\L^z,\braces{\L^{x_n}}).
\]
By Theorem \ref{circle_thm} we know that $\M_\L(\L^z,\braces{\L^{x_n}})\ge k$, establishing \eqref{AaH_cor5.4_2}.

Assume \eqref{AaH_cor5.4_2}. By \cite[Lemma~1.2]{Archbold-Kaniuth1999} there exists a sequence $\braces{\pi_n}$ converging to $\L^z$ such that $\M_\L(\L^z,\braces{\pi_n})=\M_\U(\L^z,\braces{\pi_n})=\M_\U(\L^z)$. Since the orbits are locally closed, by \cite[Proposition~5.1]{Clark2007} the mapping $G^{(0)}\rightarrow C^*(G)^\wedge:x\mapsto\L^x$ is a surjection. So there is a sequence $\braces{\L^{x_n}}$ in $C^*(G)^\wedge$ such that $\L^{x_n}$ is unitarily equivalent to $\pi_n$ for each $n$. Then
\[
\M_\L(\L^z,\braces{\L^{x_n}})\ge\M_\L(\L^z,\braces{\pi_n})=\M_\U(\L^z)\ge k,
\]
and it follows from Theorem \ref{circle_thm} that $\braces{x_n}$ is $k$-times convergent in $G^{(0)}/G$ to $z$.
\end{proof}
\end{cor}

\begin{cor}\label{AaH_cor5.7}
Suppose that $G$ is a secound-countable locally-compact Hausdorff principal groupoid with Haar system $\lambda$. Let $z\in G^{(0)}$ and let $\braces{x_n}\subset G^{(0)}$ be a sequence converging to $z$. Assume that $[z]$ is locally closed. Then the following are equivalent:
\begin{enumerate}\renewcommand{\labelenumi}{(\arabic{enumi})}
\item\label{AaH_cor5.7_1} there exists an open neighbourhood $V$ of $z$ such that $G_z^V$ is relatively compact and
\[
\underset{\scriptstyle n}{\lim\,\sup}\,\lambda_{x_n}(G^V)<\infty;
\]
\item\label{AaH_cor5.7_2} $\M_\U(\L^z,\braces{\L^{x_n}})<\infty$.
\end{enumerate}
\begin{proof}
Suppose that \eqref{AaH_cor5.7_1} holds. Since $C^*(G)$ is separable, it follows from \cite[Lemma~A.1]{Archbold-anHuef2006} that there exists a subsequence $\braces{x_{n_j}}$ of $\braces{x_n}$ such that
\[
\M_\L(\L^z,\braces{\L^{x_{n_j}}})=\M_\U(\L^z,\braces{\L^{x_{n_j}}})=\M_\U(\L^z,\braces{\L^{x_n}}).
\]
By \eqref{AaH_cor5.7_1} and Corollary \ref{AaH_cor5.6}, $\M_\L(\L^z,\braces{\L^{x_n}})<\infty$. Hence $\M_\U(\L^z,\braces{\L^{x_n}})<\infty$, as required.

Suppose that \eqref{AaH_cor5.7_1} fails. Let $\braces{V_i}$ be a basic decreasing sequence of open neighbourhoods of $z$ such that $G_z^{V_1}$ is relatively compact (such neighborhoods exist by \cite[Lemma~4.1(1)]{Clark-anHuef2010-preprint}). Then 
\[
\underset{\scriptstyle n}{\lim\,\sup}\,\lambda_{x_n}(G^{V_i})=\infty\quad\text{for each }i
\]
and we may choose a subsequence $\braces{x_{n_i}}$ of $\braces{x_n}$ such that $\lambda_{x_{n_i}}(G^{V_i})\rightarrow\infty$ as $i\rightarrow\infty$.

Let $V$ be any open neighbourhood of $z$ such that $G_z^V$ is relatively compact. There exists $i_0$ such that $V_i\subset V$ for all $i\ge i_0$. Then, for $i\ge i_0$,
\[
\lambda_{x_{n_i}}(G^{V_i})\le\lambda_{x_{n_i}}(G^V).
\]
Thus $\lambda_{x_{n_i}}(G^V)\rightarrow\infty$ as $i\rightarrow\infty$. By Corollary \ref{AaH_cor5.6}, $\M_\L(\L^z,\braces{\L^{x_n}})=\infty$. Hence $\M_\U(\L^z,\braces{\L^{x_n}})=\infty$, that is \eqref{AaH_cor5.7_2} fails.
\end{proof}
\end{cor}

\begin{cor}\label{AaH_cor5.8}
Suppose $G$ is a second-countable locally-compact Hausdorff principal groupoid with Haar system $\lambda$ such that all the orbits are locally closed. Let $z\in G^{(0)}$. Then the following are equivalent:
\begin{enumerate}\renewcommand{\labelenumi}{(\arabic{enumi})}
\item\label{AaH_cor5.8_1} $\M_\U(\L^z)<\infty$;
\item\label{AaH_cor5.8_2} there exists an open neighbourhood $V$ of $z$ such that $G_z^V$ is relatively compact and
\[
\sup_{x\in V}\lambda_x(G^V)<\infty.
\]
\end{enumerate}
\begin{proof}
If \eqref{AaH_cor5.8_2} holds then \eqref{AaH_cor5.8_1} holds by Corollary \ref{AaH_cor3.7}.

Let $\braces{V_i}$ be a basic decreasing sequence of open neighbourhoods of $z$ such that $G_z^{V_1}$ is relatively compact. If \eqref{AaH_cor5.8_2} fails then $\sup_{x\in V_i}\braces{\lambda_x(G^{V_i})}=\infty$ for each $i$ and we may choose a sequence $\braces{x_i}$ such that $x_i\in V_i$ for all $i$ and $\lambda_{x_i}(G^{V_i})\rightarrow\infty$. Since $\braces{V_i}$ is a basic decreasing sequence, $x_i\rightarrow z$.

Let $V$ be an open neighbourhood of $z$ such that $G_z^V$ is relatively compact. There exists $i_0$ such that $V_i\subset V$ for all $i\ge i_0$. Then, for $i\ge i_0$,
\[
\lambda_{x_i}(G^{V_i})\le \lambda_{x_i}(G^V).
\]
Thus $\lambda_{x_i}(G^V)\rightarrow\infty$. By Corollary \ref{AaH_cor5.7}, $\M_\U(\L^z,\braces{\L^{x_i}})=\infty$. Hence $\M_\U(\L^z)=\infty$, and so \eqref{AaH_cor5.8_1} fails.
\end{proof}
\end{cor}

\section{Graph algebra examples}\label{section_graph_algebra_examples}
We begin this section by introducing the notion of a directed graph as well as some related concepts as in the expository book \cite{Raeburn2005}, although some notation is also taken from \cite{kprr1997}. A {\em directed graph} $E=(E^0,E^1,r,s)$ consists of two countable sets $E^0$, $E^1$ and functions $r,s:E^1\rightarrow E^0$. The elements of $E^0$ and $E^1$ are called {\em vertices} and {\em edges} respectively. For each edge $e$, call $s(e)$ the {\em source} of $e$ and $r(e)$ the {\em range} of $e$. A directed graph $E$ is {\em row finite} if $r^{-1}(v)$ is finite for every $v\in E^0$.

A {\em finite path} in a directed graph $E$ is a finite sequence $\alpha=\alpha_1\alpha_2\cdots\alpha_k$ of edges $\alpha_i$ with $s(\alpha_j)=r(\alpha_{j+1})$ for $1\le j\le k-1$; write $s(\alpha)=s(\alpha_k)$ and $r(\alpha)=r(\alpha_1)$, and call $|\alpha|:=k$ the {\em length} of $\alpha$. An {\em infinite path} $x=x_1x_2\cdots$ is defined similarly, although $s(x)$ remains undefined. Let $E^*$ and $E^\infty$ denote the set of all finite paths and infinite paths in $E$ respectively. If $\alpha=\alpha_1\cdots\alpha_k$ and $\beta=\beta_1\cdots\beta_j$ are finite paths then, provided $s(\alpha)=r(\beta)$, let $\alpha\beta$ be the path $\alpha_1\cdots\alpha_k\beta_1\cdots\beta_j$. When $x\in E^\infty$ with $s(\alpha)=r(x)$ define $\alpha x$ similarly. A {\em cycle} is a finite path $\alpha$ of non-zero length such that $s(\alpha)=r(\alpha)$.

When $v$ is a vertex, $f$ is an edge, and there is exactly one infinite path with range $v$ that includes the edge $f$, then we denote this infinite path by $[v,f]^\infty$. When there is exactly one finite path $\alpha$ with $r(\alpha)=v$ and $\alpha_{|\alpha|}=f$, we denote $\alpha$ by $[v,f]^*$. In \cite{kprr1997} two paths $x,y\in E^\infty$ are defined to be {\em shift equivalent} with lag $k\in\ZZ$ (written $x\sim_k y$) if there exists $N\in\NN$ such that $x_i=y_{i+k}$ for all $i\ge N$.

Suppose $E$ is a row-finite directed graph. We refer to the groupoid constructed from $E$ by Kumjian, Pask, Raeburn and Renault in \cite{kprr1997} as the {\em path groupoid}. Before describing this construction we caution that we are using the now standard  notation for directed graphs which has the range and source swapped from the notation used in \cite{kprr1997}. This new convention is due to the development of the  higher-rank graphs, where  edges become morphisms in a category and the new convention  ensures that ``composition of morphisms is compatible with multiplication of operators in $B(\Hh)$'' \cite[p.~2]{Raeburn2005}. The path groupoid $G=G_E$ constructed from $E$ is defined as follows:
\[
G:=\braces{(x,k,y)\in E^\infty\times \ZZ\times E^\infty :x\sim_k y}.
\]
For elements of 
\[
G^{(2)}:=\braces{\big((x,k,y),(y,l,z)\big):(x,k,y),(y,l,z)\in G},
\]
Kumjian, Pask, Raeburn, and Renault defined
\[
(x,k,y)\cdot (y,l,z):=(x,k+l,z),
\]
and for arbitrary $(x,k,y)\in G$, defined
\[
(x,k,y)^{-1}:=(y,-k,x).
\]
For each $\alpha,\beta\in E^*$ with $s(\alpha)=s(\beta)$, let $Z(\alpha,\beta)$ be the set
\[
\braces{(x,k,y):x\in Z(\alpha), y\in Z(\beta), k=|\beta|-|\alpha|, x_i=y_{i+k}\text{ for } i>|\alpha|}.
\]
By \cite[Proposition~2.6]{kprr1997}, the collection of sets
\[
\braces{Z(\alpha,\beta):\alpha,\beta\in E^*, s(\alpha)=s(\beta)}
\]
is a basis of compact open sets for a second-countable locally-compact Hausdorff topology on $G$ that makes $G$ r-discrete. Kumjian, Pask, Raeburn and Renault equipped $G$ with the Haar system consisting of counting measures, which they observe is possible by first showing that a Haar system exists for the groupoid with \cite[Proposition~I.2.8]{Renault1980}, and then using \cite[Lemma~I.2.7]{Renault1980} to show that they can choose the system of counting measures.

By \cite[Corollary~2.2]{kprr1997}, the cylinder sets
\[
Z(\alpha):=\braces{x\in E^\infty :x_1=\alpha_1,\ldots,x_{|\alpha|}=\alpha_{|\alpha|}}
\]
parameterised by $\alpha\in E^*$ form a basis of compact open sets for a locally-compact $\sigma$-compact totally-disconnected Hausdorff topology on $E^\infty $. After identifying each $(x,0,x)\in G^{(0)}$ with $x\in E^\infty $,  \cite[Proposition~2.6]{kprr1997} tells us that the topology on $G^{(0)}$ is identical to the topology on $E^\infty $.

For a row-finite directed graph $E$, Kumjian, Pask, Raeburn and Renault use the path groupoid $G$ to construct the usual groupoid $C^*$-algebra $C^*(G)$, and show how a collection of partial isometries subject to some relations derived from $E$ generate $C^*(G)$. More recently,  a $C^*$-algebra $C^*(E)$ is constructed  from a collection of partial isometries subject to slightly weakened relations derived from $E$. The slightly weakened relations permit non-zero partial isometries to be related to sources in the graph, and as a result $C^*(E)$ is isomorphic to $C^*(G)$ only when $E$ contains no sources. It turns out that $C^*(E)$ and $C^*(G)$ can be substantially different: an example in \cite{Hazlewood_graph_algebra_classification} describes a graph with sources where $C^*(G)$ has continuous trace while $C^*(E)$ does not. 
In this paper we are only interested in groupoid $C^*$-algebras, so we will make no further mention of the graph algebra $C^*(E)$.

Since we wish to apply Theorem \ref{circle_thm} to path groupoids, we must be able to show that the path groupoids we consider are principal.
\begin{prop}\label{prop_principal_iff_no_cycles}
Suppose $E$ is a row-finite directed graph. The path groupoid $G$ constructed from $E$ is principal if and only if $E$ contains no cycles.
\begin{proof}
We first show that if $E$ contains no cycles then $G$ is principal. Suppose $G$ is not principal. Then there exist $x,y\in E^\infty $ and distinct $\gamma,\delta\in G$ such that $r(\gamma)=r(\delta)=x$ and $s(\gamma)=s(\delta)=y$. It follows that there exist $a,b\in\ZZ$ such that $\gamma=(x,a,y)$ and $\delta=(x,b,y)$. Notice that since $\gamma\ne\delta$, $a\ne b$. We may assume without loss of generality that $a>b$.

Now $\gamma=(x,a,y)$ implies $x\sim_a y$ and $\delta=(x,b,y)$ implies $x\sim_b y$, so there exists $N$ such that
\[
n\ge N\implies x_n=y_{n+a}=y_{n+b},
\]
and so $x_n=y_{n+a}=y_{n+a-b+b}=x_{n+a-b}$. Thus $E$ contains a cycle of length at most $a-b$.

We now show that if $G$ is principal then $E$ contains no cycles. Suppose $E$ contains the cycle $\alpha=\alpha_1\alpha_2\cdots \alpha_k$. Then $x:=\alpha\alpha\cdots$ is in $E^\infty $ with $x\sim_k x$, so both $(x,0,x)$ and $(x,k,x)$ are in $G$. It follows that $G$ is not principal.
\end{proof}
\end{prop}



\begin{example}[$2$-times convergence in a path groupoid]\label{2-times_convergence_example}
Let $E$ be the graph
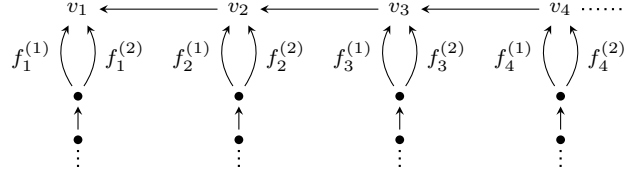
\[
\begin{alternativegraphic}[2]{merged_fun_with_groupoids_graphics}
\begin{tikzpicture}[>=stealth,baseline=(current bounding box.center)] 
\def\cellwidth{5.5};
\clip (-5em,-5.6em) rectangle (3*\cellwidth em + 4.5em,0.3em);

\foreach \x in {1,2,3,4} \foreach \y in {0} \node (x\x y\y) at (\cellwidth*\x em-\cellwidth em,-3*\y em) {$\scriptstyle v_{\x}$};
\foreach \x in {1,2,3,4} \foreach \y in {1} \node (x\x y\y) at (\cellwidth*\x em-\cellwidth em,-3*\y em) {};
\foreach \x in {1,2,3,4} \foreach \y in {2} \node (x\x y\y) at (\cellwidth*\x em-\cellwidth em,-1.5em-1.5*\y em) {};
\foreach \x in {1,2,3,4} \foreach \y in {1,2} \fill[black] (x\x y\y) circle (0.15em);


\foreach \x in {1,2,3,4} \draw [<-, bend left] (x\x y0) to node[anchor=west] {$\scriptstyle f_{\x}^{(2)}$} (x\x y1);
\foreach \x in {1,2,3,4} \draw [<-, bend right] (x\x y0) to node[anchor=east] {$\scriptstyle f_{\x}^{(1)}$} (x\x y1);

\foreach \x / \z in {1/2,2/3,3/4} \draw[black,<-] (x\x y0) to node[anchor=south] {} (x\z y0);

\foreach \x in {1,2,3,4} \draw [<-] (x\x y1) -- (x\x y2);
\foreach \x in {1,2,3,4} {
	\node (endtail\x) at (\cellwidth*\x em-\cellwidth em, -6em) {};
	\draw [dotted,thick] (x\x y2) -- (endtail\x);
}

\node(endtailone) at (3*\cellwidth em + 2.5em,0em) {};
\draw[dotted,thick] (x4y0) -- (endtailone);
\end{tikzpicture}
\end{alternativegraphic}
\]
and let $G$ be the path groupoid. For each $n\ge 1$ define $x^{(n)}:=[v_1, f_n^{(1)}]^\infty$ and let $z$ be the infinite path with range $v_1$ that passes through each $v_n$. Then $\{x^{(n)}\}$ converges $2$-times in $G^{(0)}/G$ to $z$.

\begin{proof}
We will describe two sequences in $G$ as in Definition \ref{def_k-times_convergence}. For each $n\ge 1$ define $\gamma_n^{(1)}:=(x^{(n)},0,x^{(n)})$ and $\gamma_n^{(2)}:=\big([v_1, f_n^{(2)}]^\infty,0,x^{(n)}\big)$. It follows immediately that $s(\gamma_n^{(1)})=x^{(n)}=s(\gamma_n^{(2)})$ for all $n$ and that both $r(\gamma_n^{(1)})$ and $r(\gamma_n^{(2)})$ converge to $z$ as $n\rightarrow\infty$. It remains to show that $\gamma_n^{(2)}(\gamma_n^{(1)})^{-1}\rightarrow\infty$ as $n\rightarrow\infty$.

Let $K$ be a compact subset of $G$. Our goal is to show that $\gamma_n^{(2)}(\gamma_n^{(1)})^{-1}=\gamma_n^{(2)}$ is eventually not in $K$. Since sets of the form $Z(\alpha,\beta)$ for some $\alpha,\beta\in E^*$  form a basis for the topology on the path groupoid, for each $\gamma\in K$ there exist $\alpha^{(\gamma)},\beta^{(\gamma)}\in E^*$ with $s(\alpha^{(\gamma)})=s(\beta^{(\gamma)})$ so that $Z(\alpha^{(\gamma)},\beta^{(\gamma)})$ is an open neighbourhood of $\gamma$ in $G$. Thus $\cup_{\gamma\in K}Z(\alpha^{(\gamma)},\beta^{(\gamma)})$ is an open cover of the compact set $K$, and so admits a finite subcover $\cup_{i=1}^I Z(\alpha^{(i)},\beta^{(i)})$. 

We now claim that for any fixed $n\in\PP$, if there exists $i$ with $1\le i\le I$ such that $\gamma_n^{(2)}\in Z(\alpha^{(i)},\beta^{(i)})$, then $\big|[v_1,f_n^{(2)}]^*\big|\le \abs{\alpha^{(i)}}$. Temporarily fix $n\in\PP$ and suppose there exists $i$ with $1\le i\le I$ such that $\gamma_n^{(2)}\in Z(\alpha^{(i)},\beta^{(i)})$. Suppose the converse: that $ \abs{\alpha^{(i)}}<\big|[v_1,f_n^{(2)}]^*\big|$. Since $\gamma_n^{(2)}\in Z(\alpha^{(i)},\beta^{(i)})$, it follows that $r(\gamma_n^{(2)})=[v_1, f_n^{(2)}]^\infty\in Z(\alpha^{(i)})$, and so $\alpha_p^{(i)}=[v_1, f_n^{(2)}]^\infty_p$ for every $1\le p\le|\alpha^{(i)}|$. By examining the graph we can see that $s\big([v_1, f_n^{(2)}]^\infty_p\big)=v_{p+1}$ for all $1\le p<\big|[v_1,f_n^{(2)}]^*\big|$. Since we also know that $|\alpha^{(i)}|<\big|[v_1,f_n^{(2)}]^*\big|$, we can deduce that $s(\alpha^{(i)})=v_j$ for some $j$.
Furthermore since $s(\alpha^{(i)})=s(\beta^{(i)})$, $s(\beta^{(i)})=v_j$. There is only one path with source $v_j$ and range $v_1$, so $\alpha^{(i)}=\beta^{(i)}$.
Note that when $k=|\beta^{(i)}|-|\alpha^{(i)}|$, the set $Z(\alpha^{(i)},\beta^{(i)})$ is by definition equal to
\[
\{(x,k,y):x\in Z(\alpha^{(i)}),y\in Z(\beta^{(i)}), x_p=y_{p+k} \text{ for }p>|\alpha^{(i)}|\},
\]
so since $\gamma_n^{(2)}\in Z(\alpha^{(i)},\beta^{(i)})$ and $\alpha^{(i)}=\beta^{(i)}$, we can see that $s(\gamma_n^{(2)})_p=r(\gamma_n^{(2)})_p$ for all $p>|\alpha^{(i)}|$. We know $s(\gamma_n^{(2)})=[v_1, f_n^{(1)}]^\infty$ and $r(\gamma_n^{(2)})=[v_1, f_n^{(2)}]^\infty$, so $[v_1, f_n^{(2)}]^\infty_p=[v_1, f_n^{(1)}]^\infty_p$ for all $p>|\alpha^{(i)}|$. In particular, since we assumed that $\big|[v_1,f_n^{(2)}]^*\big|>|\alpha^{(i)}|$, we have 
\[
[v_1, f_n^{(2)}]^\infty_{\big|[v_1,f_n^{(2)}]^*\big|}=[v_1, f_n^{(1)}]^\infty_{\big|[v_1,f_n^{(2)}]^*\big|},\]
so that $f_n^{(2)}=f_n^{(1)}$. But $f_n^{(1)}$ and $f_n^{(2)}$ are distinct, so we have found a contradiction, and we must have $\big|[v_1,f_n^{(2)}]^*\big|\le |\alpha^{(i)}|$.

Our next goal is to show that each $Z(\alpha^{(i)},\beta^{(i)})$ contains at most one $\gamma_n^{(2)}$. Fix $n,m\in\PP$ and suppose that both $\gamma_n^{(2)}$ and $\gamma_m^{(2)}$ are in $Z(\alpha^{(i)},\beta^{(i)})$ for some $i$. We will show that $n=m$. Since $\gamma_n^{(2)}\in Z(\alpha^{(i)}\,\beta^{(i)})$, $r(\gamma_n^{(2)})=[v_1, f_n^{(2)}]^\infty\in Z(\alpha^{(i)})$. Thus there exists $x\in E^\infty $ such that $[v_1,f_n^{(2)}]^*x\in Z(\alpha^{(i)})$ and, since $\big|[v_1,f_n^{(2)}]^*\big|\le \abs{\alpha^{(i)}}$, we can crop $x$ to form a finite $\epsilon\in E^*$ such that $[v_1,f_n^{(2)}]^*\epsilon=\alpha^{(i)}$. Similarly there exists $\delta\in E^*$ such that $[v_1,f_m^{(2)}]^*\delta=\alpha^{(i)}$. Then
\[
[v_1,f_n^{(2)}]^*\epsilon=\alpha^{(i)}=[v_1,f_m^{(2)}]^*\delta,
\]
which we can see by looking at the graph is only possible if $n=m$. We have thus shown that if $\gamma_n^{(2)}$ and $\gamma_m^{(2)}$ are in $Z(\alpha^{(i)},\beta^{(i)})$, then $\gamma_n^{(2)}=\gamma_m^{(2)}$.

Let $S=\{n\in\PP:\gamma_n^{(2)}\in K\}$. Since $K\subset \cup_{i=1}^IZ(\alpha^{(i)},\beta^{(i)})$ and since $\gamma_n^{(2)},\gamma_m^{(2)}\in Z(\alpha^{(i)},\beta^{(i)})$ implies $n=m$, $S$ can contain at most $I$ elements. Then $S$ has a maximal element $n_0$ and $\gamma_n^{(2)}\notin K$ provided $n>n_0$. Thus $\gamma_n\rightarrow\infty$ as $n\rightarrow\infty$, and we have shown that $x^{(n)}$ converges $2$-times to $z$ in $G^{(0)}/G$.
\end{proof}
\end{example}

\begin{example}[$k$-times convergence in a path groupoid]\label{k-times_convergence_example}
For any fixed positive integer $k$, let $E$ be the graph
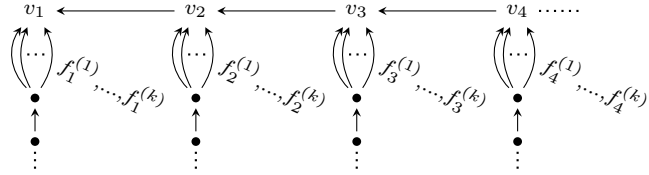
\[
\begin{alternativegraphic}[3]{merged_fun_with_groupoids_graphics}
\begin{tikzpicture}[>=stealth,baseline=(current bounding box.center)] 
\def\cellwidth{5.5};
\clip (-5em,-5.6em) rectangle (3*\cellwidth em + 4.5em,0.3em);

\foreach \x in {1,2,3,4} \foreach \y in {0} \node (x\x y\y) at (\cellwidth*\x em-\cellwidth em,-3*\y em) {$\scriptstyle v_{\x}$};
\foreach \x in {1,2,3,4} \foreach \y in {1} \node (x\x y\y) at (\cellwidth*\x em-\cellwidth em,-3*\y em) {};
\foreach \x in {1,2,3,4} \foreach \y in {2} \node (x\x y\y) at (\cellwidth*\x em-\cellwidth em,-1.5em-1.5*\y em) {};
\foreach \x in {1,2,3,4} \foreach \y in {1,2} \fill[black] (x\x y\y) circle (0.15em);


\foreach \x in {1,2,3,4} \draw [<-, bend right=40] (x\x y0) to node[anchor=west] {} (x\x y1);
\foreach \x in {1,2,3,4} \draw [<-, bend right=25] (x\x y0) to node[anchor=west] {} (x\x y1);

\foreach \x in {1,2,3,4} \draw [<-, bend left] (x\x y0) to node[anchor=west,rotate=-25] {$\scriptstyle f_{\x}^{(1)},\ldots,f_{\x}^{(k)}$} (x\x y1);
\foreach \x in {1,2,3,4} \draw [dotted,thick] ($(x\x y0)!{1/(2*cos(10))}!-10:(x\x y1)$) -- ($(x\x y0)!{1/(2*cos(14))}!14:(x\x y1)$);

\foreach \x / \z in {1/2,2/3,3/4} \draw[black,<-] (x\x y0) to node[anchor=south] {} (x\z y0);

\foreach \x in {1,2,3,4} \draw [<-] (x\x y1) -- (x\x y2);
\foreach \x in {1,2,3,4} {
	\node (endtail\x) at (\cellwidth*\x em-\cellwidth em, -6em) {};
	\draw [dotted,thick] (x\x y2) -- (endtail\x);
}

\node(endtailone) at (3*\cellwidth em + 2.5em,0em) {};
\draw[dotted,thick] (x4y0) -- (endtailone);
\end{tikzpicture}
\end{alternativegraphic}
\]
and let $G$ be the path groupoid. For each $n\ge 1$ define $x^{(n)}:=[v_1, f_n^{(1)}]^\infty$ and let $z$ be the infinite path that passes through each $v_n$. Then the sequence $\{x^{(n)}\}$ converges $k$-times in $G^{(0)}/G$ to $z$.
\begin{proof}
After defining $\gamma_n^{(i)}:=\big([v_1, f_n^{(i)}]^\infty,0,x^{(n)}\big)$ for each $1\le i\le k$, an argument similar to that in Example \ref{2-times_convergence_example} establishes the $k$-times convergence.
\end{proof}
\end{example}

\begin{example}\label{ML2_MU3_example}[Lower multiplicity 2 and upper multiplicity 3]
Consider the graph $E$ described by
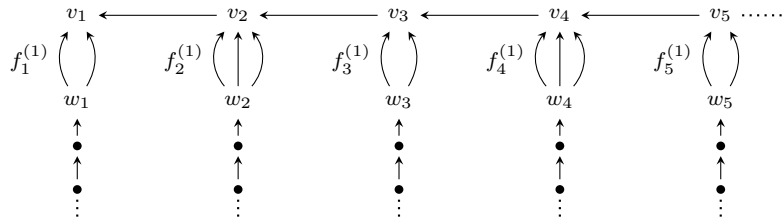
\[
\begin{alternativegraphic}[3]{merged_fun_with_groupoids_graphics}
\begin{tikzpicture}[>=stealth,baseline=(current bounding box.center)] 
\def\cellwidth{5.5};
\clip (-2.5em,-7.1em) rectangle (4*\cellwidth em + 2.5em,0.3em);

\foreach \x in {1,2,3,4,5} \foreach \y in {0} \node (x\x y\y) at (\cellwidth*\x em-\cellwidth em,-3*\y em) {$\scriptstyle v_{\x}$};
\foreach \x in {1,2,3,4,5} \foreach \y in {1} \node (x\x y\y) at (\cellwidth*\x em-\cellwidth em,-3*\y em) {$\scriptstyle w_{\x}$};
\foreach \x in {1,2,3,4,5} \foreach \y in {2,3} \node (x\x y\y) at (\cellwidth*\x em-\cellwidth em,-1.5em-1.5*\y em) {};
\foreach \x in {1,2,3,4,5} \foreach \y in {2,3} \fill[black] (x\x y\y) circle (0.15em);

\foreach \x in {1,3,5} \draw [<-, bend left] (x\x y0) to node[anchor=west] {} (x\x y1);
\foreach \x in {1,3,5} \draw [<-, bend right] (x\x y0) to node[anchor=east] {$\scriptstyle f_{\x}^{(1)}$} (x\x y1);
\foreach \x in {2,4} \draw [<-, bend left=40] (x\x y0) to node[anchor=west] {} (x\x y1);
\foreach \x in {2,4} \draw [<-, bend right=40] (x\x y0) to node[anchor=east] {$\scriptstyle f_{\x}^{(1)}$} (x\x y1);
\foreach \x in {2,4} \draw [<-] (x\x y0) to node {} (x\x y1);

\foreach \x / \z in {1/2,2/3,3/4,4/5} \draw[black,<-] (x\x y0) to node[anchor=south] {} (x\z y0);

\foreach \x in {1,2,3,4,5} \draw [<-] (x\x y1) -- (x\x y2);
\foreach \x in {1,2,3,4,5} \draw [<-] (x\x y2) -- (x\x y3);
\foreach \x in {1,2,3,4,5} {
	\node (endtail\x) at (\cellwidth*\x em-\cellwidth em, -7.5em) {};
	\draw [dotted,thick] (x\x y3) -- (endtail\x);
}

\node(endtailone) at (4*\cellwidth em + 2.5em,0em) {};
\draw[dotted,thick] (x5y0) -- (endtailone);
\end{tikzpicture}
\end{alternativegraphic}
\]
where for each odd $n\ge 1$ there are exactly two paths $f_n^{(1)},f_n^{(2)}$ with source $w_n$ and range $v_n$, and for each even $n\ge 2$ there are exactly three paths $f_n^{(1)},f_n^{(2)},f_n^{(3)}$ with source $w_n$ and range $v_n$. Let $G$ be the path groupoid, define $x^{(n)}:=[v_1, f_n^{(1)}]^\infty$ for every $n\ge 1$, and let $z$ be the infinite path that meets every vertex $v_n$ (so $z$ has range $v_1$). Then
\[
\M_\L(\L^z,\braces{\L^{x^{(n)}}})=2\quad\text{and}\quad\M_\U(\L^z,\braces{\L^{x^{(n)}}})=3.
\]
\begin{proof}
We know that $\braces{x^{(n)}}$ converges $2$-times to $z$ in $G^{(0)}/G$ by the argument in Example \ref{2-times_convergence_example}, so we can apply Theorem \ref{circle_thm} to see that $\M_\L(\L^z,\braces{\L^{x^{(n)}}})\ge 2$. We can see that the subsequence $\braces{x_{2n}}$ of $\braces{x^{(n)}}$ converges $3$-times to $z$ in $G^{(0)}/G$ by Example \ref{k-times_convergence_example}. Theorem \ref{2nd_circle_thm} now tells us that $\M_\U(\L^z,\braces{\L^{x^{(n)}}})\ge 3$.

Now suppose $\M_\L(\L^z,\braces{\L^{x^{(n)}}})\ge 3$. Then by Theorem \ref{circle_thm}, $\braces{x^{(n)}}$ converges $3$-times to $z$ in $G^{(0)}/G$, so there must exist three sequences $\braces{\gamma_n^{(1)}}$,$\braces{\gamma_n^{(2)}}$, and $\braces{\gamma_n^{(3)}}$ as in the definition of $k$-times convergence (Definition \ref{def_k-times_convergence}). For each odd $n$, there are only two elements in $G$ with source $x^{(n)}$, so there must exist $1\le i<j\le 3$ such that $\gamma_n^{(i)}=\gamma_n^{(j)}$ frequently. Then $\gamma_n^{(j)}(\gamma_n^{(i)})^{-1}=r(\gamma_n^{(i)})$ frequently and, since $r(\gamma_n^{(i)})\rightarrow z$, $\braces{\gamma_n^{(j)}(\gamma_n^{(i)})^{-1}}$ admits a convergent subsequence. Thus $\gamma_n^{(j)}(\gamma_n^{(i)})^{-1}\nrightarrow\infty$, contradicting the definition of $k$-times convergence.

If $\M_\U(\L^z,\braces{\L^{x^{(n)}}})\ge 4$, then by Theorem \ref{2nd_circle_thm} there is a subsequence of $\braces{x^{(n)}}$ that converges $4$-times to $z$ in $G^{(0)}/G$. A similar argument to that in the preceding paragraph shows that this is not possible since there are at most 3 edges between any $v_n$ and $v_m$. It follows that $\M_\L(\L^z,\braces{\L^{x^{(n)}}})=2$ and $\M_\U(\L^z,\braces{\L^{x^{(n)}}})=3$.
\end{proof}
\end{example}

\begin{lemma}\label{upper_and_lower_multiplicities_lemma}
In Example \ref{2-times_convergence_example},
\[
\M_\L(\L^z,\braces{\L^{x^{(n)}}})=\M_\U(\L^z,\braces{\L^{x^{(n)}}})=2;
\]
and in Example \ref{k-times_convergence_example},
\[
\M_\L(\L^z,\braces{\L^{x^{(n)}}})=\M_\U(\L^z,\braces{\L^{x^{(n)}}})=k.
\]
\begin{proof}
The same argument as that found in Example \ref{ML2_MU3_example} can be used to demonstrate this lemma. The explicit proof was given for Example \ref{ML2_MU3_example} since it covers the case where the upper and lower multiplicities are distinct.
\end{proof}
\end{lemma}

In the next example we will add some structure to the graph from Example \ref{2-times_convergence_example} to create a path groupoid $G$ with non-Hausdorff orbit space that continues to exhibit $2$-times convergence.  
\begin{example}\label{example_with_non-Hausdorff_orbit_space}
Let $E$ be the directed graph
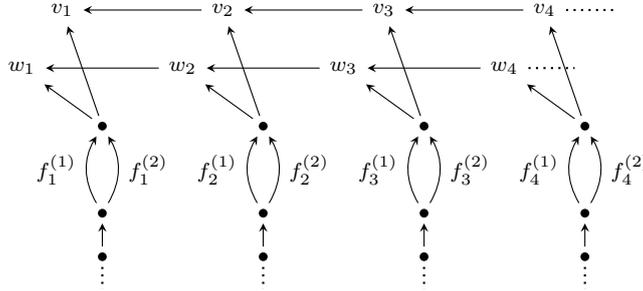
\[
\begin{tikzpicture}[>=stealth,baseline=(current bounding box.center)] 
\def\cellwidth{5.5};
\clip (-4em,-5.6em) rectangle (3*\cellwidth em + 3em,4.3em);

\foreach \x in {1,2,3,4} \foreach \y in {0,1} \node (x\x y\y) at (\cellwidth*\x em-\cellwidth em,-3*\y em) {};
\foreach \x in {1,2,3,4} \foreach \y in {2} \node (x\x y\y) at (\cellwidth*\x em-\cellwidth em,-1.5em-1.5*\y em) {};
\foreach \x in {1,2,3,4} \foreach \y in {0,1,2} \fill[black] (x\x y\y) circle (0.15em);

\foreach \x in {1,2,3,4} {
	\node (x\x z1) at (\cellwidth*\x em-1.25*\cellwidth em,4em) {$\scriptstyle v_\x$};
	\node (x\x z2) at (\cellwidth*\x em-1.5*\cellwidth em,2em) {$\scriptstyle w_\x$};
}
\foreach \x in {1,2,3,4} {
	\draw[<-] (x\x z1) -- (x\x y0);
	\draw[<-] (x\x z2) -- (x\x y0);
}

\foreach \x in {1,2,3,4} \draw [<-, bend left] (x\x y0) to node[anchor=west] {$\scriptstyle f_{\x}^{(2)}$} (x\x y1);
\foreach \x in {1,2,3,4} \draw [<-, bend right] (x\x y0) to node[anchor=east] {$\scriptstyle f_{\x}^{(1)}$} (x\x y1);

\foreach \x / \z in {1/2,2/3,3/4} \draw[black,<-] (x\x z1) to node[anchor=south] {} (x\z z1);
\foreach \x / \z in {1/2,2/3,3/4} \draw[black,<-] (x\x z2) to node[anchor=south] {} (x\z z2);

\foreach \x in {1,2,3,4} \draw [<-] (x\x y1) -- (x\x y2);
\foreach \x in {1,2,3,4} {
	\node (endtail\x) at (\cellwidth*\x em-\cellwidth em, -6em) {};
	\draw [dotted,thick] (x\x y2) -- (endtail\x);
}

\node (endtailone) at ($(x4z1) + 0.5*(\cellwidth em, 0 em)$) {};
\node (endtailtwo) at ($(x4z2) + 0.5*(\cellwidth em, 0 em)$) {};
\draw[dotted,thick] (x4z1) -- (endtailone);
\draw[dotted,thick] (x4z2) -- (endtailtwo);
\end{tikzpicture}
\]
and let $G$ be the path groupoid. For every $n\ge 1$ let $x^{(n)}$ be the infinite path $[v_1, f_n^{(1)}]^\infty$. Let $x$ be the infinite path with range $v_1$ that passes through each $v_n$ and let $y$ be the infinite path with range $w_1$ that passes through each $w_n$. Then the orbit space $G^{(0)}/G$ is not Hausdorff and $\braces{x^{(n)}}$ converges $2$-times in $G^{(0)}/G$ to both $x$ and $y$.

\begin{proof}
To see that $\braces{x^{(n)}}$ converges $2$-times to $x$ in $G^{(0)}/G$, consider the sequences $\braces{([v_1, f_n^{(2)}]^\infty,0,x^{(n)})}$ and $\braces{(x^{(n)},0,x^{(n)})}$ and follow the argument as in Example \ref{2-times_convergence_example}. To see that $\braces{x^{(n)}}$ converges $2$-times to $y$ in $G^{(0)}/G$, consider the sequences $\braces{([w_1,f_n^{(1)}]^\infty,0,x^{(n)})}$ and $\braces{([w_1,f_n^{(2)}]^\infty,0,x^{(n)})}$. While it is tempting to think that this example exhibits $4$-times convergence (or even $3$-times convergence), this is not the case (see Example~\ref{ML2_MU3_example} for an argument demonstrating this). We know $x^{(n)}$ converges $k$-times to $x$ in $G^{(0)}/G$, so $[x^{(n)}]\rightarrow [x]$ in $G^{(0)}/G$, and similarly $[x^{(n)}]\rightarrow [y]$ in $G^{(0)}/G$. It follows that $G^{(0)}/G$ is not Hausdorff since $[x]\ne [y]$.
\end{proof}
\end{example}

In all of the examples above, the orbits in  $G^{(0)}$ are closed and hence $C^*(G)$ and $G^{(0)}/G$ are homeomorphic by \cite[Proposition~5.1]{Clark2007}.  By combining the features of the graphs in Examples~\ref{ML2_MU3_example} and~\ref{example_with_non-Hausdorff_orbit_space} we obtain a principal groupoid whose $C^*$-algebra has non-Hausdorff spectrum and  distinct upper and lower multiplicities among its irreducible representations.


\begin{thebibliography}{10}

\bibitem{Archbold1994}
Robert~J. Archbold.
\newblock Upper and lower multiplicity for irreducible representations of
  {$C^\ast$}-algebras.
\newblock {\em Proc. London Math. Soc.}, 69(1):121--143, 1994.


\bibitem{Archbold-Deicke2005}
Robert~J. Archbold and Klaus Deicke.
\newblock Bounded trace {$C\sp *$}-algebras and integrable actions.
\newblock {\em Math. Z.}, 250(2):393--410, 2005.

\bibitem{Archbold-anHuef2006}
Robert~J. Archbold and Astrid an~Huef.
\newblock Strength of convergence in the orbit space of a transformation group.
\newblock {\em J. Funct. Anal.}, 235(1):90--121, 2006.

\bibitem{Archbold-anHuef2008}
Robert~J. Archbold and Astrid an~Huef.
\newblock Strength of convergence and multiplicities in the spectrum of a
  {$C^*$}-dynamical system.
\newblock {\em Proc. London Math. Soc.}, 96(3):545--581, 2008.


\bibitem{Archbold-Kaniuth1999}
Robert~J. Archbold and Eberhard Kaniuth.
\newblock Upper and lower multiplicity for irreducible representations of
  {SIN}-groups.
\newblock {\em Illinois J. Math.}, 43(4):692--706, 1999.

\bibitem{aklss2001}
Robert~J. Archbold, Eberhard Kaniuth, Jean Ludwig, G\"{u}nter Schlichting, and
  Douglas W.~B. Somerset.
\newblock Strength of convergence in duals of {$C^*$}-algebras and nilpotent
  {L}ie groups.
\newblock {\em Adv. Math.}, 158(1):26--65, 2001.

\bibitem{Archbold-Ludwig-Schlichting2007}
Robert~J. Archbold, Jean Ludwig, and G\"{u}nter Schlichting.
\newblock Limit sets and strengths of convergence for sequences in the duals of
  thread-like {L}ie groups.
\newblock {\em Math. Z.}, 255(2):245--282, 2007.

\bibitem{Archbold-Somerset-Spielberg1997}
Robert~J. Archbold, Douglas W.~B. Somerset, and Jack~S. Spielberg.
\newblock Upper multiplicity and bounded trace ideals in {$C\sp *$}-algebras.
\newblock {\em J. Funct. Anal.}, 146(2):430--463, 1997.

\bibitem{Archbold-Spielberg1996}
Robert~J. Archbold and Jack~S. Spielberg.
\newblock Upper and lower multiplicity for irreducible representations of
  {$C\sp *$}-algebras. {II}.
\newblock {\em J. Operator Theory}, 36(2):201--231, 1996.

\bibitem{Clark2007}
Lisa~O. Clark.
\newblock Classifying the types of principal groupoid {$C\sp *$}-algebras.
\newblock {\em J. Operator Theory}, 57(2):251--266, 2007.


\bibitem{Clark-anHuef2008}
Lisa~O. Clark and Astrid an~Huef.
\newblock Principal groupoid {$C\sp *$}-algebras with bounded trace.
\newblock {\em Proc. Amer. Math. Soc.}, 136(2):623--634, 2008.


\bibitem{Clark-anHuef2010-preprint}
Lisa~O. Clark and Astrid an~Huef.
\newblock The representation theory of $C^*$-algebras associated to groupoids.
\newblock Preprint, 2010.


\bibitem{Green1977}
Philip Green.
\newblock {$C\sp*$}-algebras of transformation groups with smooth orbit space.
\newblock {\em Pacific J. Math.}, 72(1):71--97, 1977.

\bibitem{Hazlewood_graph_algebra_classification}
Robert Hazlewood.
\newblock Continuous trace, {F}ell, bounded trace, liminal and postliminal
  graph algebras.
\newblock In preparation.

\bibitem{kprr1997}
Alex Kumjian, David Pask, Iain Raeburn, and Jean Renault.
\newblock Graphs, groupoids, and {C}untz-{K}rieger algebras.
\newblock {\em J. Funct. Anal.}, 144(2):505--541, 1997.

\bibitem{Ludwig1990}
Jean Ludwig.
\newblock On the behaviour of sequences in the dual of a nilpotent {L}ie group.
\newblock {\em Math. Ann.}, 287(2):239--257, 1990.

\bibitem{Muhly-book}
Paul~S. Muhly.
\newblock Coordinates in operator algebra.
\newblock (Book in preparation).

\bibitem{Muhly-Renault-Williams1996}
Paul~S. Muhly, Jean~N. Renault, and Dana~P. Williams.
\newblock Continuous-trace groupoid {$C^\ast$}-algebras. {III}.
\newblock {\em Trans. Amer. Math. Soc.}, 348(9):3621--3641, 1996.

\bibitem{Muhly-Williams1990}
Paul~S. Muhly and Dana~P. Williams.
\newblock Continuous trace groupoid {$C\sp *$}-algebras.
\newblock {\em Math. Scand.}, 66(2):231--241, 1990.

\bibitem{Pedersen1989}
Gert~K. Pedersen.
\newblock {\em Analysis now}, volume 118 of {\em Graduate Texts in
  Mathematics}.
\newblock Springer-Verlag, New York, 1989.

\bibitem{Raeburn2005}
Iain Raeburn.
\newblock {\em Graph algebras}, volume 103 of {\em CBMS Regional Conference
  Series in Mathematics}.
\newblock Published for the Conference Board of the Mathematical Sciences,
  Washington, DC, 2005.

\bibitem{Ramsay1990}
Arlan Ramsay.
\newblock The {M}ackey-{G}limm dichotomy for foliations and other {P}olish
  groupoids.
\newblock {\em J. Funct. Anal.}, 94(2):358--374, 1990.

\bibitem{Renault1980}
Jean Renault.
\newblock {\em A groupoid approach to {$C\sp{\ast} $}-algebras}, volume 793 of
  {\em Lecture Notes in Mathematics}.
\newblock Springer, Berlin, 1980.

\bibitem{Williams2007}
Dana~P. Williams.
\newblock {\em Crossed products of {$C{^\ast}$}-algebras}, volume 134 of {\em
  Mathematical Surveys and Monographs}.
\newblock American Mathematical Society, Providence, RI, 2007.

\end{thebibliography}
\end{document}